\documentclass[12pt]{article}

\usepackage[letterpaper,hmargin=27mm,vmargin=23mm]{geometry}
\usepackage{amsfonts,amsmath,amsthm,amssymb,bm,bbm,dsfont}
\usepackage[utf8]{inputenc}
\usepackage{lmodern}
\usepackage{upgreek}
\usepackage{enumitem}
\usepackage{natbib,color}
\usepackage[pdftex]{graphicx}
\usepackage{subfigure,float,booktabs}
\usepackage{bbm}
\usepackage{tikz-cd}
\usepackage{quiver}
\usetikzlibrary{arrows,automata}
\usepackage{threeparttable}
\definecolor{winered}{rgb}{0.5,0,0}
\usepackage[pdfstartview=FitH,  bookmarks=true,  bookmarksopen=true,  breaklinks=true,  colorlinks=true,
linkcolor=winered,anchorcolor=my,  citecolor=blue, filecolor=blue,  menucolor=blue,  urlcolor=black]{hyperref} 
\allowdisplaybreaks

\numberwithin{equation}{section}
\newtheorem{theorem}{Theorem}[section]

\newtheorem{proposition}[theorem]{Proposition}
\newtheorem{corollary}[theorem]{Corollary}

\newtheorem{remark}[theorem]{Remark}
\theoremstyle{definition}

{\theoremstyle{plain}
	}
\definecolor{my}{rgb}{0.05,0.05,0.5}
\definecolor{myBlue}{rgb}{.1,.1,.5}
\definecolor{myGreen}{rgb}{0,.4,0}
\definecolor{myRed}{rgb}{.25,0.15,.5}

\definecolor{my}{rgb}{0.05,0.05,0.5}

\newcommand{\cond}{\displaystyle \stackrel{d}{\longrightarrow}}

\newcommand{\conas}{\stackrel{a.s.}{\longrightarrow}}
\newcommand{\conp}{\stackrel{p}{\longrightarrow}}

\makeatletter
\renewcommand\paragraph{%
	\@startsection{paragraph}{4}{0mm}%
	{-\baselineskip}%
	{.2\baselineskip}%
	{\normalfont\normalsize\bfseries}}
\makeatother

\renewcommand{\liminf}{\displaystyle \operatornamewithlimits{\lim\inf \ }}
\renewcommand{\mathbf}[1]{\textbf{\textit{#1}}}

\newcommand{\Norm}[1]{\mathcal{N}\left(#1\right)}

\newcommand{\rs}{\mathcal{S}}

\newcommand{\E}{\operatorname{E}}
\newcommand{\V}{\operatorname{Var}}

\newcommand{\Rmnum}[1]{\expandafter\@slowromancap\romannumeral #1@}

\begin{document}
		{\title{{A Combinatorial Central Limit Theorem for Stratified Randomization}\thanks{The author thanks Xavier D'Haultf\oe{}uille for helpful discussions that led to this work.}}
		\date{}
		\author{
			Purevdorj Tuvaandorj\thanks{York University, {tpujee@yorku.ca.}}
		}
	}
	\maketitle
\begin{abstract}
This paper establishes a combinatorial central limit theorem for stratified randomization, which holds under a Lindeberg-type condition. The theorem allows for an arbitrary number or sizes of strata, with the sole requirement being that each stratum contains at least two units. This flexibility accommodates both a growing number of large and small strata simultaneously, while imposing minimal conditions. We then apply this result to derive the asymptotic distributions of two test statistics proposed for instrumental variables settings in the presence of potentially many strata of unrestricted sizes.
\begin{description}
\item[Keywords:] Combinatorial central limit theorem, permutation test, randomization inference, stratification.
\end{description}
\end{abstract}
\section{Introduction}
Stratified randomization is frequently utilized to enhance covariate balance in randomized experiments. 
In this scheme, blocks or strata are initially created based on pre-treatment variables, and randomization inference is subsequently conducted by randomizing units within each stratum \citep[see e.g.,][Chapter~9]{Imbens-Rubin(2015)}. Stratified randomization can be used for obtaining exact tests \citep{Imbens-Rosenbaum(2005), Andrews-Marmer(2008), Zhao-Ding2021, DHT2023} and efficient estimators \citep{Cytrynbaum2021, Armstrong2022, BLSTM2023}, and is a natural inferential choice for models with survey data obtained through stratified sampling.\par 

The key tool used when establishing the asymptotic validity of the permutation tests is 
the combinatorial central limit theorem (CLT) \citep{Hoeffding1951}. 
While the standard combinatorial CLT or finite population CLT can be applied to derive the asymptotic distributions of stratified permutation statistics when there is a finite number of large strata
\citep{Imbens-Rosenbaum(2005), LiDing2017}, its application to general stratified permutations is not straightforward. This complication arises because the number of strata, either large or small, or both, may \emph{grow} with the sample size (for instance, consider the case of having $n^{1/2}$ many strata, each of size $n^{1/2}$, where $n$ denotes the sample size).\par 

Despite recent contributions by \cite{LiuYang2020} and \cite{LiuRenYang2022} that propose finite population CLTs tailored for such scenarios, a combinatorial CLT accommodating potentially many strata of unrestricted size remains essential for addressing general stratified permutation/rank statistics. In this context, this paper introduces a general form of the combinatorial CLT for stratified randomization, applicable under a \emph{Lindeberg-type} condition. The proof is based on Stein's method \citep[see, e.g.,][for a detailed account]{CGS2011}, and in particular, draws from the constructions of \cite{Bolthausen1984} and \cite{Schneller1988}. Leveraging this result, we characterize the asymptotic distribution of a test statistic proposed by \cite{Imbens-Rosenbaum(2005)} and a variant of a test statistic by \cite{Andrews-Marmer(2008)} in the instrumental variables (IVs) setting featuring many strata, under weak assumptions.\par 

A substantial body of literature exists on stratified randomized experiments and large sample results for stratified randomization. For empirical practices and an overview, see, for example, \cite{BruhnMcKenzie2009}, \cite{Imbens-Rubin(2015)} and \cite{Ding2023}. Papers that consider a coarse stratification involving a finite number of strata, where asymptotic normality is typically derived via finite population or rank CLTs, include \cite{Andrews-Marmer(2008)}, \cite{LiDing2017}, \cite{Zhao-Ding2021}. Papers that consider a finely stratified experiment with many strata, where the asymptotic normality is derived via Lindeberg or Lyapunov CLTs, include \cite{Fogarty2018mit}, \cite{Bai2022} and \cite{BaiRomanoShaikh2022}.\par
 
Furthermore, while \cite{LiuYang2020} and \cite{LiuRenYang2022} provide finite population CLTs for stratified experiments involving many strata, \cite{DHT2023} establish a combinatorial CLT with many strata to examine the properties of a stratified permutation subvector test in linear regressions. The CLT for stratified randomization in this paper is sufficiently general to encompass the aforementioned CLTs as it handles both large and small strata simultaneously and requires a weaker Lindeberg-type condition. In particular, finite population CLTs with many strata follow from the result akin to the derivation of usual finite population CLTs from standard combinatorial or rank CLTs.\par 

Stratified permutation is an instance of permutations with \emph{restricted positions} \citep[see, e.g.,][]{Rosenbaum1984,Diaconis2001,LeiBickel2021}; thus, the set of such permutations forms a subset of all permutations of observation indices. \citet[][Chapter 6]{CGS2011} provide normal approximation results for a class of restricted permutations {\textemdash} in contrast with the set of entire permutations required in \citet{Hoeffding1951}'s CLT {\textemdash} but do not consider stratified permutations. The present paper fills this gap in the literature by developing normal approximation for this important class of permutations.\par

This paper proceeds as follows. Section \ref{sec: Main} lays out the basic framework and presents the combinatorial CLT for stratified randomization with an arbitrary number of strata. Section \ref{sec: link FPCLT} establishes a connection between our result and finite population CLTs. Section \ref{sec: IV} applies the combinatorial CLT to derive the asymptotic distribution of two test statistics in IV settings. The proofs are provided in the appendix.
\section{Main result}\label{sec: Main}
Let $n$ be a fixed sample size, ${S}$ be an integer-valued random variable (${S}\ge 1$ a.s.), and $s=1,\dots,{S}$ denote strata of sizes $n_s\ge 2$ a.s., with $\sum_{s=1}^{{S}} n_s=n$. Consider a double array of real-valued random variables in stratum $s$: $\{a_{ij}^s\}_{i,j=1}^{n_s}$, and define $\bm{a}_{n}\equiv \{a_{ij}^s\}_{s=1,\dots, {S}; (i,j)\in\{1,\dots, n_s\}^2}$. Here, the distributions of $\{a_{ij}^s\}_{i,j=1}^{n_s}$ and ${S}$ may depend on $n$, and the stratification may be based on random auxiliary covariates. \par 
Let $\pi$ be a \emph{stratified} permutation that permutes the indices within each stratum: 

\begin{equation*}
\pi=
\underbrace{\begin{pmatrix}
	1&\dots &n_1\\
	\pi_1(1)&\dots &\pi_1(n_1)
\end{pmatrix}}_{\text{Stratum $1$}}
\dots
\underbrace{
\begin{pmatrix}
	1&\dots &n_{{S}}\\
	\pi_S(1)&\dots &\pi_{{S}}(n_{{S}})
\end{pmatrix}}_{\text{Stratum ${S}$}},
\end{equation*}
where $\pi_s$ is a permutation of $\{1,\dots, n_s\}$ in stratum $s=1,\dots,{S}$.\footnote{For simplicity, we do not further index the elements of $\{1,\dots, n_s\}$ by $s$ throughout the paper.} 
We denote by $\mathbb{S}_n$ the set of all such permutations with  $|\mathbb{S}_n|=\prod_{s=1}^{{S}}n_s!$. In this paper, we consider permutations uniformly distributed over $\mathbb{S}_n$ for $\bm{a}_{n}$ given: $\pi\sim\mathcal{U}(\mathbb{S}_n)$, where $\mathcal{U}(A)$ denotes the uniform distribution on a finite set $A$. Let $P^\pi$ be the probability measure of $\pi$ conditional on $\bm{a}_{n}$, and $\E_\pi[\cdot]$ and $\V_\pi[\cdot]$ denote the corresponding expectation and variance operators.

Furthermore, we use the notation $1(\cdot)$ for the indicator function, and  
$\Phi(\cdot)$ and $\Phi^{-1}(\cdot)$ for the cumulative distribution and quantile function of a standard real  normal distribution. $\sum_{i,j=1}^{n_s}$ and $\sum_{i\neq j}$ are the shortcuts for $\sum_{i=1}^{n_s}\sum_{j=1}^{n_s}$ and $\sum_{i,j=1; i\neq j}^{n_s}$, respectively. We abbreviate ``left-hand side'' and ``right-hand side'' as LHS and RHS. Let $\Vert\cdot\Vert$ denote the Frobenius norm, $\lambda_{\min}(\cdot)$ denote the minimum eigenvalue of a square matrix, and $\max_{s, i}$ denote for the maximum over $s\in\{1,\dots,{S}\}$ and $i\in\{1,\dots,n_s\}$ for a given $n$. \par 
\par 

The main result of this paper is the following combinatorial CLT for stratified randomization which is based on a Lindeberg-type condition.
\begin{theorem}[Combinatorial CLT for stratified randomization]\label{HoeffdingCLT}
	Let ${S}$ be an integer-valued random variable such that ${S}\ge 1$ a.s. and $s=1,\dots,{S}$ denote strata of sizes $n_s\ge 2$ a.s., with $\sum_{s=1}^{{S}} n_s=n$. Let $\pi\sim\mathcal{U}(\mathbb{S}_n)$ and for each $s$, $\{a_{ij}^s\}_{i,j=1}^{n_s}$ be $n_s^2$ scalar random variables satisfying:
	\begin{enumerate}[label=(\alph*)]
		\item\label{cond: centered}$\sum_{i=1}^{n_s}a_{ij}^{s}=\sum_{j=1}^{n_s}a_{ij}^s=0$ for $1\leq i,j\leq n_s$ and  $s=1,\dots, {S}$ a.s.;
		\item\label{cond: Lindeberg} For any $\epsilon>0$, $\sigma_n^{-2}\sum_{s=1}^{{S}} \frac{1}{n_s}\left(\sum_{i,j=1}^{n_s}a_{ij}^{s2}1(\sigma_n^{-1}|a_{ij}^s|>\epsilon)\right)\conas 0$ as $n\to\infty$, where 
		$\sigma_n^2\equiv \sum_{s=1}^{{S}} \frac{1}{n_s-1}\left(\sum_{i,j=1}^{n_s}a_{ij}^{s2}\right)$ and $a_{ij}^s\neq 0$ for some $(i,j)$ and $s$. 
	\end{enumerate}
	Let $T^\pi\equiv \sum_{s=1}^{{S}} \sum_{i=1}^{n_s} a^{s}_{i\pi(i)}/\sigma_n$. Then, as $n\to\infty$
	\begin{equation}\label{eq: CCLT}
	P^\pi(T^\pi\leq t)\conas \Phi(t)\ \ \text{for any}\ \ t\in\mathbb{R}.
	\end{equation}
\end{theorem}
\begin{remark}\label{remark: R}~
\normalfont
\begin{enumerate}[label={(\arabic{enumi})},leftmargin=*]
\item\label{R1}When ${S}=1$ with probability 1, $n_{{S}}=n$ and the result specializes to Theorem 1 of \cite{Motoo1956}. The latter provides a version of \cite{Hoeffding1951}'s CLT based on a Lindeberg-type condition that corresponds to the condition in \ref{cond: Lindeberg}.
\item\label{R2} The key ingredient in the proof of Theorem \ref{HoeffdingCLT} is the random selection of a stratum $s$ with probability $p_s\equiv n_s/n$. We then bound $|\E_\pi[T^\pi f(T^\pi)]-\E_\pi[f'(T^\pi)]|$, where $f(\cdot)$ is defined in \eqref{def: f} below. This is achieved by perturbing the sum  $\sum_{i=1}^{n_s} a^{s}_{i\pi(i)}$ within the stratum $s$, following the approach in \cite{Bolthausen1984} and \cite{Schneller1988}, and then ``averaging'' over the strata $s=1,\dots, S$.  
\item\label{R3} Since $n_s\geq 2$ a.s., ${S}\leq n/2$ and $n-{S}\geq n/2$ a.s. On noting that $n-{S}=\sum_{s=1}^{S}(n_s-1)\leq \prod_{s=1}^{{S}}n_s\leq \vert\mathbb{S}_n\vert$ a.s. and $n/2\to \infty$, we have $|\mathbb{S}_n|\conas \infty$.  
\item\label{R4} A Lyapunov-type sufficient condition for the Lindeberg condition in \ref{cond: Lindeberg} is 
	\begin{equation}\label{cond: Lyap}
\frac{1}{\sigma_n^{2+\delta}}\sum_{s=1}^{{S}} \frac{1}{n_s}\sum_{i,j=1}^{n_s}|a_{ij}^{s}|^{2+\delta}\conas 0 
	\end{equation} 
as $n\to\infty$, for some $\delta>0$.
\item\label{R5} The randomness of $T^\pi$ stems from $\bm{a}_{n}$ and $\pi$. If $\bm{a}_{n}$ are random variables such that the condition \ref{cond: centered} holds a.s. and the convergence in \ref{cond: Lindeberg} holds in probability, then using a subsequencing argument, the convergence in \eqref{eq: CCLT} can be restated as 
\begin{equation}\label{eq: CCLT2}
P^\pi(T^\pi\leq t)\conp \Phi(t). 
\end{equation}
\item\label{R6} 
Consider the case $a_{ij}^s=\tilde{b}_{si}\tilde{c}_{sj}$, where $\tilde{b}_{si}\equiv b_{si}-\bar{b}_s$, $\bar{b}_s\equiv n_{s}^{-1}\sum_{i=1}^{n_s}b_{si}$, $\tilde{c}_{sj}\equiv c_{sj}-\bar{c}_s$, $\bar{c}_s\equiv n_{s}^{-1}\sum_{j=1}^{n_s}c_{sj}$ and $b_{si}$ and $c_{sj}$ are scalars.
\cite{DHT2023} (see Lemma 4 therein) establish a combinatorial CLT for stratified randomization based on the following assumptions:  as $n\to\infty$
\begin{enumerate}[label=(\alph*)]
	\item\label{cond: sigcon} For $\sigma_n^2\equiv\sum_{s=1}^{{S}} \frac{1}{n_s-1}\left(\sum_{i=1}^{n_s}\tilde{b}_{si}^{2}\right)\left(\sum_{i=1}^{n_s}\tilde{c}_{si}^{2}\right)$,
	$\sigma_n^2\conp \sigma^2>0$;
	\item \label{cond: 3mom} $\sum_{s=1}^{{S}}n_s^{-1} \left(\sum_{i=1}^{n_s} |\tilde{b}_{si}|^3\right)
	\left(\sum_{i=1}^{n_s} |\tilde{c}_{si}|^3\right)\conp 0$ ;
	\item\label{cond: Uvarcon}
	$\sum_{s=1}^{{S}}(n_s-1)^{-1}\left(\sum_{i=1}^{n_s}\tilde{b}_{si}^{4}\right)\left(\sum_{i=1}^{n_s}\tilde{c}_{si}^{4}\right)\conp 0$.
\end{enumerate}
From the conditions \ref{cond: sigcon} and \ref{cond: 3mom}, we have 
\begin{equation}\label{cond: Lyap2}
	\frac{1}{\sigma_n^{3}}\sum_{s=1}^{{S}} \frac{1}{n_s}\sum_{i,j=1}^{n_s}|a_{ij}^{s}|^{3}
	=\frac{1}{\sigma_n^{3}}\sum_{s=1}^{{S}} \frac{1}{n_s}\left(\sum_{i=1}^{n_s}|\tilde{b}_{si}|^{3}\right)\left(\sum_{i=1}^{n_s}|\tilde{c}_{si}|^{3}\right)\conp 0.	
\end{equation} 
Thus, the condition in \eqref{cond: Lyap} holds in probability with $\delta=1$, and \eqref{eq: CCLT2} follows from a subsequencing argument.
\end{enumerate}
\end{remark}
\noindent A leading special case of interest similar to that considered in Remark \ref{remark: R}\ref{R6} is the asymptotic normality of the sum $n^{-1/2}\sum_{s=1}^{{S}}\sum_{i=1}^{n_s}\tilde{b}_{si}\tilde{c}_{s\pi(i)}$, where $\tilde{b}_{si}\equiv b_{si}-\bar{b}_s$, $\bar{b}_s\equiv n_{s}^{-1}\sum_{i=1}^{n_s}b_{si}$, $b_{si}\in\mathbb{R}^k$ with $k\geq 1$,  $\tilde{c}_{si}\equiv c_{si}-\bar{c}_s$, $c_{si}\in\mathbb{R}$ and $\bar{c}_s\equiv n_{s}^{-1}\sum_{i=1}^{n_s}c_{si}$. The following result holds as a corollary to Theorem \ref{HoeffdingCLT}. 

\begin{corollary}\label{cor: vector bc}
Suppose that $\lambda_{\min}(\Sigma_n)>\lambda>0$ for $\Sigma_n\equiv n^{-1}\sum_{s=1}^{{S}} \frac{1}{n_s-1}\left(\sum_{i=1}^{n_s}\tilde{b}_{si}\tilde{b}_{si}'\right)\left(\sum_{j=1}^{n_s}\tilde{c}_{sj}^{2}\right)$,  and for some positive constants $\delta$ and $M_0$, either
	\begin{enumerate}[label=(\alph*)]
	\item\label{cond: Vector bc1}
$n^{-1}\sum_{s=1}^{{S}}\sum_{i=1}^{n_s}\Vert {b}_{si}\Vert^{4+\delta}<M_0$ a.s. and 
$n^{-1}\sum_{s=1}^{{S}}\sum_{i=1}^{n_s}|{c}_{si}|^{4+\delta}<\infty$ a.s. for all $n$; or
\item\label{cond: Vector bc2} $n^{-1/2}\max_{s,i}\Vert b_{si}\Vert\conas 0$, 
$n^{-1}\sum_{s=1}^{{S}}\sum_{i=1}^{n_s}\Vert {b}_{si}\Vert^{2}<M_0$ a.s. 
and $\max_{1\leq s\leq S}n_s^{-1}\sum_{i=1}^{n_s}\vert {c}_{si}\vert^{2+\delta}<M_0$ a.s. for all $n$.
\end{enumerate}
Then, as $n\to\infty$ 
\begin{equation}\label{eq: Vector AN} 
\Sigma_n^{-1/2}n^{-1/2}\sum_{s=1}^{{S}}\sum_{i=1}^{n_s}\tilde{b}_{si}\tilde{c}_{s\pi(i)}\cond \Norm{0, I_k}\quad\text{a.s}.
\end{equation}
\end{corollary}
\section{Relation to finite population CLT for stratified randomization}\label{sec: link FPCLT}
This section explores the connection between the combinatorial CLT of Theorem \ref{HoeffdingCLT} and a finite population CLT. 
\cite{Bickel1984} and \cite{LiuYang2020} establish finite population CLTs for linear combinations of stratum means under stratified sampling, where the number of strata may diverge.\par 
Consider a finite population divided into non-random ${S}$ strata: $(y_{s1},\dots, y_{sn_s}), s=1,\dots, {S}$. The stratum-specific mean and variance are defined respectively as  $\bar{y}_s\equiv n_s^{-1}\sum_{i=1}^{n_s}y_{si}$ and $v_{s}^2\equiv \frac{1}{n_s-1}\sum_{i=1}^{n_s}(y_{si}-\bar{y}_s)^2$, while the population mean and weighted variance are given by 
\begin{equation}\label{eq: popmeanvar}
\bar{y}\equiv n^{-1}\sum_{s=1}^S\sum_{i=1}^{n_s}y_{si}=\sum_{s=1}^Sp_s\bar{y}_s\ \ \text{and}\ \
v^2\equiv \sum_{s=1}^{S}p_s^2v_s^2\frac{n_s-n_{1s}}{n_{1s}n_s},
\end{equation}
where $p_s\equiv n_s/n$ and $n_{1s}, 2\leq n_{1s}\leq n_s-1,$ denotes the number of units sampled without replacement from stratum $s$ in \cite{Bickel1984}'s set-up (the 
number of treated units in stratum $s$ in \cite{LiuYang2020}'s set-up). The sampling indicators are  $(Z_{s1},\dots, Z_{sn_s})\in\{0,1\}^{n_s}, s=1,\dots, S$, where $Z_{si}=1$ if the unit $i$ in stratum $s$ is sampled, and $0$ otherwise. The probability that $\{(Z_{s1},\dots,Z_{sn_s})\}_{s=1}^S$ takes on a value
$\{(z_{s1},\dots,z_{sn_s})\}_{s=1}^S$ such that $n_{1s}=\sum_{i=1}^{n_s}z_{si}$ for $s=1,\dots, S,$  is $\prod_{s=1}^Sn_{1s}!(n_s-n_{1s})!/n_s!$. The weighted sample mean and variance are then defined as
\begin{equation*}
 \hat{y}\equiv \sum_{s=1}^Sp_s\hat{y}_s\ \ \text{and}\ \ \hat{v}^2\equiv \sum_{s=1}^Sp_s^2\hat{v}_{s}^2\frac{n_s-n_{1s}}{n_{1s}n_s},
\end{equation*}
where $\hat{y}_s\equiv n_{1s}^{-1}\sum_{i=1}^{n_s}Z_{si}y_{si}$ and $\hat{v}_{s}^2
\equiv \frac{1}{n_s-1}\sum_{i=1: Z_{si}=1}^{n_s}(y_{si}-\hat{y}_s)^2$ are the stratum-specific sample mean and variance.
With the notations introduced, the Lindeberg condition employed by \cite{LiuYang2020s} (Condition A1 and Theorem A1 therein) may be restated as follows: for any $\epsilon>0$
\begin{equation}\label{eq: LYLind}
	\sum_{s=1}^{{S}} \frac{1}{n_s-1}\left(\sum_{i=1}^{n_s}w_s^2\frac{(y_{si}-\bar{y}_s)^2}{v_s^2}1\left(w_s\frac{|y_{si}-\bar{y}_s|}{v_s}>\epsilon \sqrt{\frac{n_{1s}(n_s-n_{1s})}{n_s}}\right)\right)\to 0
\end{equation}
as $n\to\infty$, where 
\begin{equation}\label{eq: ws Vs}
	w_s^2\equiv \frac{p_s^2v_s^2(n_s-n_{1s})}{v^2n_sn_{1s}}.
\end{equation}
Under the condition in \eqref{eq: LYLind}, \cite{Bickel1984} show that $(\hat{y}-\bar{y})/{\hat{v}}\cond \Norm{0,1}$ as $n\to\infty$. To draw a link between the latter result and Theorem \ref{HoeffdingCLT}, we let $a_{ij}^s=\tilde{b}_{si}\tilde{c}_{sj}$ as in 
Remark \ref{remark: R}\ref{R6} and Corollary \ref{cor: vector bc}, where $b_{si}$ and $c_{sj}$ are now defined as 
\begin{equation}\label{eq: LYb}
	b_{si}\equiv 
	\begin{cases}
		\frac{1}{n_{1s}}&\text{if}\quad i=1,\dots, n_{1s},\\
		0&\text{if}\quad i=n_{1s}+1,\dots, n_s
	\end{cases},
	\quad c_{sj}\equiv p_sy_{sj}.
\end{equation}
The choice of $b_{si}$ in \eqref{eq: LYb} appears in \cite{Madow1948} and \cite{Sen1995}, among others. 
By simple algebra, $\bar{b}_s=n_s^{-1}\sum_{i=1}^{n_s}b_{si}=n_s^{-1}$, $\bar{c}_s= n_s^{-1}\sum_{j=1}^{n_s}c_{sj}=p_s\bar{y}_s$, and 
\begin{equation}\label{eq: LYb2}
	\tilde{b}_{si}= b_{si}-\bar{b}_s=
	\begin{cases}
		\frac{1}{n_{1s}}-\frac{1}{n_s}&\text{if}\quad i=1,\dots, n_{1s},\\
		-\frac{1}{n_s}&\text{if}\quad i=n_{1s}+1,\dots, n_s
	\end{cases},\quad 
	\tilde{c}_{sj}= c_{sj}-\bar{c}_s= p_s({y}_{sj}-\bar{y}_s).
\end{equation}
It turns out that the variance
$\sigma_n^2$ defined in the condition \ref{cond: Lindeberg} of Theorem \ref{HoeffdingCLT} and $v^2$ in \eqref{eq: popmeanvar} are identical:
\begin{equation}\label{eq: vareq}
\sigma_n^2
=\sum_{s=1}^{{S}} \frac{1}{n_s-1}\left(\sum_{i=1}^{n_s}\tilde{b}_{si}^{2}\right)\left(\sum_{i=1}^{n_s}\tilde{c}_{si}^{2}\right)=\sum_{s=1}^{{S}} \frac{1}{n_s-1}\frac{(n_s-n_{1s})p_s^2}{n_{1s}n_s}\left(\sum_{i=1}^{n_s}({y}_{si}-\bar{y}_s)^{2}\right)=v^2.
\end{equation}
Then, through a direct calculation given in Section \ref{sec:R2}, the Lindeberg condition in Theorem  \ref{HoeffdingCLT}\ref{cond: Lindeberg} follows from the condition in \eqref{eq: LYLind}, as summarized in the remark below. 
\begin{remark}\label{remark: R2}
\normalfont
The condition in \eqref{eq: LYLind} implies the Lindeberg condition in Theorem \ref{HoeffdingCLT}\ref{cond: Lindeberg}. 
\end{remark}

\section{Stratified randomization inference with IVs}\label{sec: IV}
This section presents two applications of the CLT in Theorem \ref{HoeffdingCLT}. First, we establish the asymptotic distribution of a statistic proposed by 
\cite{Imbens-Rosenbaum(2005)}. Second, we show that the test of \cite{Andrews-Marmer(2008)} may be conservative in the presence of small strata. We then consider a simple modification of their statistic that yields a correct asymptotic level and derive its null asymptotic distribution.\par 

\subsection{\cite{Imbens-Rosenbaum(2005)} statistic}
\cite{Imbens-Rosenbaum(2005)} consider an IV setting with non-random $S$ strata $s=1,\dots, S$ each of size $n_s$, given by: $$Y-\beta D=r_C,$$ where  
$D=[D_1',\dots, D_{S}']\in\mathbb{R}^n$ with $D_s=[D_{s1},\dots, D_{sn_s}]'\in\mathbb{R}^{n_s}$ and $n=\sum_{s=1}^Sn_s$, denotes the vector of doses received by the units, 
$Y=[Y_1',\dots, Y_S']'\in\mathbb{R}^n$ with $Y_s=[Y_{s1},\dots, Y_{sn_s}]'\in\mathbb{R}^{n_s}$ is the vector of responses exhibited by the units, 
$r_C=[r_{C1}',\dots, r_{CS}']'\in\mathbb{R}^n$ with $r_{Cs}=[r_{Cs1},\dots, r_{Csn_s}]'\in\mathbb{R}^{n_s}$ is the vector of responses if the units received the control dose, and the scalar parameter $\beta$, with a true value $\beta_0$, measures of the dose-response relationship. 
\par 
Let $h_s=[h_{s1},\dots, h_{sn_s}]'\in\mathbb{R}^{n_s}$ be a sorted and fixed IV vector such that 
$h_{sj}\leq h_{s(j+1)}$ for each $j=1,\dots, n_s-1$ and $s=1,\dots,S$. For $h\equiv [h_1',\dots, h_S']'\in\mathbb{R}^{n}$, assume that the observed IV $Z$ is randomly assigned i.e.
$Z=h_\pi\in\mathbb{R}^n$, where $\pi\sim \mathcal{U}(\mathbb{S}_n)$ and $\mathbb{S}_n$ is the set of all stratified permutations in this setting.\par 
Here, $r_C=Y-\beta_0 D$ is assumed to be fixed, hence independent of $Z$. For testing the hypothesis $H_0:\beta=\beta_0$,  
\cite{Imbens-Rosenbaum(2005)} propose the following statistic:
\begin{equation}\label{def: T}
	T=q(Y-\beta_0D)'\rho(Z),
\end{equation}
where $q(\cdot)\in\mathbb{R}^n$ is some method of scoring responses such as the ranks or the aligned ranks, and 
$\rho(\cdot)\in\mathbb{R}^n$ is some method of scoring the instruments that satisfies $\rho_\pi(h)=\rho(h_\pi)=\rho(Z)$. Let $\rho(h)=[\rho_{1}',\dots, \rho_S']'$ with $\rho_s=[\rho_{s1},\dots, \rho_{sn_s}]'$, and $q(Y-\beta_0D)=[q_{1}',\dots, q_{S}']'$ with $q_s=[q_{s1},\dots, q_{sn_s}]'$. \par 
Before we state the result for the asymptotic distribution of the statistic $T$, we remark that 
$T=\sum_{s=1}^S\sum_{i=1}^{n_s}q_{si}\rho_{s\pi(i)}$ using the notations above. With a convention that $\V_\pi[\sum_{i=1}^{n_s}q_{si}\rho_{s\pi(i)}]=0$ for stratum of size $n_s=1$, it follows from Proposition 1 of \cite{Imbens-Rosenbaum(2005)} that 
\begin{align}
	\mu
	&\equiv \E_\pi[T]=\sum_{s=1}^S{n_s}\bar{q}_{s}\bar{\rho}_{s},\label{eq: MeanT}\\
	\sigma^2\equiv \V_\pi[T]
	&=\sum_{s=1, n_s\geq 2}^S\frac{1}{n_s-1}\left(\sum_{i=1}^{n_s}\tilde{q}_{si}^2\right)\left(\sum_{i=1}^{n_s}\tilde{\rho}_{si}^2\right),\label{eq: varT}
\end{align}
where $\bar{q}_s\equiv n_s^{-1}\sum_{i=1}^{n_s}q_{si}$, $\bar{\rho}_s\equiv n_s^{-1}\sum_{i=1}^{n_s}\rho_{si}$, 
and $\tilde{q}_{si}\equiv q_{si}-\bar{q}_s$ and $\tilde{\rho}_{si}\equiv \rho_{si}-\bar{\rho}_s$. \par 

 \cite{Imbens-Rosenbaum(2005)} outline strategies for deriving the asymptotic distribution of the statistic $T$ in the presence of either a few large strata whose sizes tend to infinity, or many small strata, whose sizes remain bounded, and suggest applying the rank-central limit theorem \citep[see e.g.,][Chapter 6]{Hajek-Sidak-Sen(1999)} to the former and standard CLTs to the latter.
 
 \cite{LiDing2017} consider the case where $q(\cdot)$ is an identity function and $\rho(Z)$ is a vector of binary IVs, and derive the null asymptotic distribution of the statistic in \eqref{def: T}. The general CLT provided in Theorem \ref{HoeffdingCLT} allows us to derive the asymptotic distribution of the statistic $T$ in \eqref{def: T} at once under general conditions, allowing for many strata without categorizing them into large and small ones or restricting their sizes. The following proposition gives a formal justification for this argument.
\begin{proposition}\label{prop: SRIVa}
Suppose that $\liminf_{n\to\infty}\sigma^2>\lambda>0$, and for some positive constants $M_0$ and $\delta$, either 
	\begin{enumerate}[label=(\alph*)]
		\item\label{Asbound2} $n^{-1}\sum_{s=1}^S\sum_{i=1}^{n_s}|\rho_{si}|^{4+\delta}<M_0$ and $n^{-1}\sum_{s=1}^S\sum_{i=1}^{n_s}|q_{si}|^{4+\delta}<M_0$ for all $n$; or 
		\item\label{Asbound1} $n^{-1/2}\max_{s,i}|\rho_{si}|\to 0$,  $n^{-1}\sum_{s=1}^S\sum_{i=1}^{n_s}|\rho_{si}|^{2}<M_0$ and $\max_{1\leq s\leq S}n_s^{-1}\sum_{i=1}^{n_s}|{q}_{si}|^{2+\delta}<M_0$ for all $n$. 
\end{enumerate} 
Then, under $H_0:\beta=\beta_0$, ${\sigma}^{-1}(T-\mu)\cond \Norm{0,1}$ as $(n-S)\to \infty$.  
\end{proposition}
\begin{remark}\label{rem: ImbensRosenbaum}~
\normalfont
\begin{enumerate}[label={(\arabic{enumi})},leftmargin=*]
		\item The strata of size $n_s=1$ are discarded in the statistic ${\sigma}^{-1}(T-\mu)$. Since $n-S$, which corresponds to the effective sample size, tends to infinity, so does the number of observations stratified into strata of sizes $2$ or greater.
		Consequently, the condition $n\to\infty$ in Theorem \ref{HoeffdingCLT} can be replaced by $(n-S)\to \infty$. 
		
		For $S$ random, a sufficient condition for $(n-S)\conas \infty$ is provided in Lemma 2 of \cite{DHT2023}. 
		\item Assumption \ref{Asbound1} on $q_{si}$ holds, for example, if $q_{si}=\varphi(i/(n_s+1))$, where $\varphi:(0,1)\mapsto \mathbb{R}$ is a score function 
		satisfying $\int_{0}^{1}\vert {\varphi}(x)\vert^{2+\delta}dx<M_0/2$. This is because $(n_s+1)/n_s\leq 2$ and 
		\begin{equation}\label{eq: int ineq}
		(n_s+1)^{-1}\sum_{i=1}^{n_s}|\varphi(i/(n_s+1))|^{2+\delta}\leq \int_{0}^{1}|\varphi(x)|^{2+\delta}dx,
		\end{equation}
		 as the LHS of the inequality above is the sum of the areas of $n_s$ 
		 rectangles with sides $|\varphi(i/(n_s+1))|^{2+\delta}$ and $[i/(n_s+1), (i+1)/(n_s+1)]$, $i=1,\dots, n_s,$ which lie between the function $|\varphi(x)|^{2+\delta}$ and $[0,1]$.\par 
		For the Wilcoxon score $\varphi(x)=x$, letting $\delta=1$, the RHS of \eqref{eq: int ineq} is  $\int_{0}^{1}x^3dx=1/4$.  For the normal (or van der Waerden) score $\varphi(x)=\Phi^{-1}(x)$, letting $\delta=2$,  we have, through a change of variables,
		$\int_{0}^{1}(\Phi^{-1}(x))^4dx=3$. So the condition on $q_{si}$ in Assumption \ref{Asbound1} is satisfied for these score functions.
	\end{enumerate}
\end{remark}
\subsection{\cite{Andrews-Marmer(2008)} statistic}
Next, we consider a related rank-based Anderson-Rubin-type statistic \citep{Anderson-Rubin(1949)} studied in \citet[][Section~3]{Andrews-Marmer(2008)} under a super population framework. 
 The model is 
 \begin{equation}\label{eq: AMmodel}
 	Y_{si}=\alpha_s+D_{si}\beta+r_{Csi},\quad i=1,\dots, n_s;\ s=1,\dots, S,
 \end{equation}
 where $\alpha_s$ denotes a stratum-specific scalar parameter, $Y_{si}$ is the response variable, $D_{si}$ is a scalar endogenous regressor, $r_{Csi}$ is the error term which, unlike the previous case, is random, and the number of strata $S$ is not random.  In addition, there is a $k$-vector of IVs, $Z_{si}$, with $k\geq 1$, that does not include
 a constant. The model in \eqref{eq: AMmodel} may arise, for example, due to stratification based on empirical support points 
 $\{X_1,\dots, X_S\}$ of auxiliary regressors $X\in\mathbb{R}^p, p\geq 1,$
  i.e. 
 $\alpha_s\equiv X_s'\eta_s$, where $X_s\in\mathbb{R}^p$ and a parameter vector $\eta_s\in\mathbb{R}^p$ are constant within stratum $s$, but may vary across strata. The hypothesis of interest is again $H_0:\beta=\beta_0$.\par 
 Let $R_{si}$ denote the rank of $Y_{si}-\beta_0 D_{si}$ among $Y_{s1}-\beta_0 D_{s1},\dots, Y_{sn_s}-\beta_0 D_{sn_s}$, and 
 $\varphi:(0,1)\mapsto \mathbb{R}$ be a score function. The 
 \cite{Andrews-Marmer(2008)} statistic (Equation (3.3) therein) may be defined as 
 \begin{equation}\label{def: AMRAR}
B_{n}
\equiv n\,A_{n}'\Omega_{n}^{-1}A_{n},
\end{equation}
where $A_{n}=\sum_{s=1}^SA_{ns}$, $A_{ns}\equiv n^{-1}\sum_{i=1}^{n_s}(Z_{si}-\bar{Z}_s)\varphi\left(\frac{R_{si}}{n_s+1}\right)$, 
$\bar{Z}_s\equiv n_s^{-1}\sum_{i=1}^{n_s}Z_{si}$, and 
\begin{equation}\label{eq: AMcovm}
	\Omega_{n}
	\equiv n^{-1}\sum_{s=1}^{S}\sum_{i=1}^{n_s}(Z_{si}-\bar{Z}_s)(Z_{si}-\bar{Z}_s)'\int_{0}^{1}(\varphi(x)-\bar{\varphi})^2dx,
\end{equation}
with $\bar{\varphi}\equiv\int_{0}^{1}\varphi(x)dx$. Remark here that if $\{r_{Csi}\}_{i=1}^{n_s}$ are i.i.d. with a continuous distribution, $\{R_{si}\}_{i=1}^{n_s}$ are uniformly distributed over $\{1,2,\dots, n_s\}$ \citep[see e.g.,][Lemma 13.1]{vanderVaart(1998)}, so we can reformulate $A_n$ as a stratified permutation statistic. We deviate from \cite{Andrews-Marmer(2008)} in a minor way by considering the statistic 
\begin{equation}\label{def: Bnstar}
B_{n}^{*}\equiv n\,A_n'\Omega_n^{*-1}A_n,
\end{equation}
with a covariance matrix estimator $\Omega_n^{*}$ defined as 
\begin{equation}\label{eq: AMmcovm}
	\Omega_{n}^{*}
	\equiv n^{-1}\sum_{s=1, n_s\geq 2}^{S}\left\{\sum_{i=1}^{n_s}(Z_{si}-\bar{Z}_s)(Z_{si}-\bar{Z}_s)'\right\}
	\left\{\frac{1}{n_s-1}\sum_{i=1}^{n_s}\left(\varphi_{si}-\bar{\varphi}_s\right)^2\right\},
\end{equation}
where ${\varphi}_{si}\equiv \varphi\left(\frac{i}{n_s+1}\right)$ and $\bar{\varphi}_s\equiv n_{s}^{-1}\sum_{i=1}^{n_s}\varphi_{si}$. 
The estimator in \eqref{eq: AMmcovm} is preferred to that in \eqref{eq: AMcovm} because 
the former is the variance of $n^{1/2}A_{n}$ (with respect to the distribution of the ranks), whereas the covariance matrix in \eqref{eq: AMcovm} is an approximation valid only when the strata are large, a scenario considered by \cite{Andrews-Marmer(2008)}. In fact, when there are many small strata, the term $\int_{0}^{1}(\varphi(x)-\bar{\varphi})^2dx$
will involve a nonnegligible error as an approximation of $(n_s-1)^{-1}\sum_{i=1}^{n_s}\left(\varphi_{si}-\bar{\varphi}_s\right)^2$. By way of example, let $n_s=2$ for all $s=1,\dots, S$ and $\varphi(x)=\Phi^{-1}(x)$. Simple calculations yield  $\bar{\varphi}_s=(\Phi^{-1}(1/3)+\Phi^{-1}(2/3))/2=0$ and 
\begin{equation}\label{eq: step}
\frac{1}{n_s-1}\sum_{i=1}^{n_s}\left(\varphi_{si}-\bar{\varphi}_s\right)^2
=(\Phi^{-1}(1/3))^2+(\Phi^{-1}(2/3))^2\approx 0.37.
\end{equation}
This is far below the value $\int_{0}^{1}(\varphi(x)-\bar{\varphi})^2dx=1$, so the covariance matrix estimator in \eqref{eq: AMcovm} is biased upward as an estimator of the asymptotic variance of $n^{1/2}A_{n}$ and the level-$\alpha$ test that rejects when $B_n$ exceeds the $1-\alpha$ quantile of $\chi^2_k$ distribution will be conservative. The LHS of \eqref{eq: step} (after rounding) is 0.56, 0.69, 0.89, 0.96, 0.98 for $n_s=5, 10, 50, 200, 500$, respectively, indicating a reduction in the approximation error with increasing strata sizes. A similar result holds for the Wilcoxon score as well.

The null asymptotic distribution of the statistic $B_{n}^{*}$ is provided in the proposition below. 

\begin{proposition}\label{prop: AMstat}
Assume that $\liminf_{n\to\infty}\lambda_{\min}(\Omega_{n}^{*})>\lambda>0$, and 
	\begin{enumerate}[label=(\alph*)]
	\item\label{AM1} $n^{-1/2}\max_{s,i}\Vert Z_{si}\Vert\to 0$, $n^{-1}\sum_{s=1}^S\sum_{i=1}^{n_s}\Vert Z_{si}\Vert^{2}<C_0$ for all $n$, and 
	$\int_{0}^{1}\vert {\varphi}(x)\vert^{2+\delta}dx<C_0$ 
	for some positive constants $\delta$ and $C_0$;
	\item\label{AM2} $\{r_{Csi}: 1\leq i\leq n_s, 1\leq s\leq S\}$ are independent random variables and for each $s$,  $\{r_{Csi}\}_{i=1}^{n_s}$ have an identical continuous distribution. 
\end{enumerate} 
Then, under $H_0:\beta=\beta_0$, $B_{n}^{*}\cond \chi^2_k$ as $(n-S)\to \infty$.  
\end{proposition}
\begin{remark}\label{rem: AndrewMarmer}~
	\normalfont
The assumptions on the IVs and error terms in \ref{AM1}  and \ref{AM2} are comparable to Assumptions C1 and 4 of \cite{Andrews-Marmer(2008)}, and the integrability assumption on the score function 
in \ref{AM1} is slightly stronger than that in Assumption 3 of \cite{Andrews-Marmer(2008)} but is not overly restrictive. 
\end{remark}
As argued above, the test based on the statistic in \eqref{def: AMRAR}, which employs a $\chi^2_k$ critical value, may be conservative in the presence of many small strata, whereas the test based on the modified statistic in \eqref{def: Bnstar} maintains an asymptotically correct level. This discrepancy is illustrated in a simple numerical experiment
based on the following design:
\begin{align*}
Y_i&=D_i\beta+\gamma_1+X_i\gamma_2+u_i,\\
D_i&=Z_i'\pi+\psi_1+X_i\psi_2+v_i,\quad i=1,\dots, n,
\end{align*}
where we first make $n$ independent draws of $Z_i\sim t_5(0, I_k)$ with $k=1,3$ corresponding to just and over-identified cases.\footnote{Here, $t_5(0, I_k)$ stands for the $k$-variate $t$-distribution with degrees of freedom $5$, zero mean and covariance matrix $I_k$.} For each $k$, we draw the stratification variable $X_i\sim \mathcal{U}(\{1,\dots, \frac{n}{r}\})/( \frac{n}{r})$, 
where $r$ varies over $\{2,5,10,25,80\}$. A smaller value of $r$ leads to a greater number of strata, $S$, equal to the empirical support size of $X_i$. We keep $\{(Z_i',X_i)'\}_{i=1}^n$ thus generated fixed over the simulation replications. For each replication, we let $u_i\sim t_1$, $\epsilon_i\sim t_1$, and $v_i=\rho u_i+\sqrt{1-\rho^2}\epsilon_i,$ where $\rho=0.5$ and $u_i$ and $\epsilon_i$ are mutually independent and i.i.d. across $i=1,\dots, n$. The values $\beta,\gamma_1,\gamma_2, \psi_1$ and $\psi_2$ are set to $0$, and $\pi=(1,\dots, 1)'\sqrt{\lambda/(nk)}\in\mathbb{R}^k$ with $\lambda=9$, indicating moderate identification strength of the IVs. The sample size is set relatively large at $n=400$ to elucidate the effect of the strata sizes embodied in $r$.\par
After stratifying the synthetic data according to the realized values of $X_i$, we compute the following four statistics to test the null hypothesis $H_0:\beta=0$: the statistic in \eqref{def: AMRAR} using both the normal and Wilcoxon score functions, denoted as $B_n^N$ and $B_n^W$, respectively, and their modified counterparts $B_n^{N*}$ and $B_n^{W*}$ based on \eqref{def: Bnstar}.\par 
Table \ref{tab: size} displays the empirical size of the tests based on 5000 replications. The maximum stratum size $\max n_s$ and the number of strata $S$ are proportional and inversely proportional to $r$ respectively.
The tests $B_n^N$ and $B_n^W$ under-reject when the strata sizes are small, although the latter shows less size distortion. Not surprisingly, the rejection rates of the tests improve as the strata size increases.  
In contrast, the modified versions $B_n^{N*}$ and $B_n^{W*}$ exhibit reasonably accurate rejection rates in all cases considered, consistent with the theoretical predictions.\footnote{The result remains qualitatively similar when $Z_i$ and $X_i$ are randomly drawn in each replication, and is available upon request.}

\begin{table}[htbp!]
	\begin{center}
		\begin{threeparttable}
			\caption{Null rejection rates for $H_0:\beta=0$ at $5\%$ level}		\label{tab: size}
		\begin{tabular}{crrrrrrrrrr}
	\toprule
&\multicolumn{5}{c}{$k=1$}&\multicolumn{5}{c}{$k=3$}\\
\cmidrule(lr){2-6} \cmidrule(lr){7-11} 
$r$ & 2 & 5 & 10 & 25 & 80 & 2 & 5 & 10 & 25 & 80 \\ 
  \midrule
$B_n^{N*}$ & 4.74 & 4.72 & 5.30 & 4.80 & 4.64 & 4.96 & 4.76 & 4.68 & 5.02 & 4.62 \\ 
$B_n^{W*}$ & 4.76 & 4.66 & 5.44 & 4.92 & 4.84 & 4.78 & 4.98 & 4.46 & 4.90 & 4.74 \\ 
$B_n^{N}$& 0.44 & 0.92 & 2.06 & 3.08 & 4.10 & 0.02 & 0.32 & 0.94 & 2.52 & 3.64 \\ 
$B_n^{W}$ & 2.12 & 2.96 & 4.20 & 4.60 & 4.76 & 1.54 & 2.32 & 3.02 & 4.16 & 4.60 \\ \midrule 
$\max n_s$ & 7.00 & 10.00 & 15.00 & 31.00 & 83.00 & 7.00 & 10.00 & 16.00 & 31.00 & 86.00 \\ 
$S$ & 162.00 & 79.00 & 40.00 & 16.00 & 5.00 & 174.00 & 80.00 & 40.00 & 16.00 & 5.00 \\ 
   \bottomrule
\end{tabular}
			\begin{tablenotes}
				\footnotesize{\item Notes: $B^N_n$ and $B_n^W$ denote the normal and Wilcoxon score rank statistics, respectively, and $B^{N*}_n$ and $B^{W*}_n$ denote their modified versions. 5000 replications.}	
			\end{tablenotes}
		\end{threeparttable}
	\end{center}
\end{table}
\newpage 
\appendix
\section{Proof of Theorem \ref{HoeffdingCLT}} 

\label{sub:lemma_ref_hoeffding_clt}
We will show that $\E_\pi[h(T^\pi)]\conas \E[h(Z)]$ as $n\to\infty$, 
where $Z\sim \Norm{0,1}$ and $h$ is a fixed continuous function that can be extended continuously to $\bar{\mathbb{R}}=\mathbb{R}\cup\{\pm \infty\}$. Letting $\phi(x)$ denote the standard normal density function, we define 
\begin{equation}\label{def: f}
	f(x)\equiv \phi(x)^{-1}\int_{-\infty}^x (h(y)-\E[h(Z)])dy.
\end{equation}
It can then be shown that \citep[see][]{Bolthausen1984}
\begin{align}
	f'(x)-xf(x)
	&=h(x)-\E[h(Z)]\quad \forall x\in\mathbb{R},\label{eq: f prop1}\\
	\lim_{x\to\pm\infty}f'(x)
	&=0\ \text{and $f'(x)$ is uniformly continuous}.\label{eq: f prop2}
\end{align}
The proof is divided into two steps.
In Step 1, we introduce perturbations of the permutation $\pi\sim\mathcal{U}(\mathbb{S}_n)$. In Step 2, we derive the asymptotic normality. 

\paragraph*{Step 1: Perturbations of $\pi$}
Let 
$\mathbb{N}_s\equiv \{1,\dots, n_s\}$ be the set of indices for observations in stratum $s=1,\dots, {S}$, and $\mathbb{G}_s$ be the set of all permutations of $\mathbb{N}_s$. 
We can then write $\pi$ as $\pi=(\pi_1,\dots,\pi_{{S}})$, where $\pi_s\sim \mathcal{U}(\mathbb{G}_s), s=1,\dots, {S}$. 

Let $\rs$ be a random variable taking values in $\{1,\dots, {S}\}$ independent of $\pi$ with $P(\rs=s)=p_s\equiv n_s/n$. In the stratum $\rs$, we choose a quadruple of indices $(I_1^{\rs}, I_2^{\rs}, J_1^{\rs}, J_2^{\rs})$ as follows: 
let $(I_1^{\rs}, I_2^{\rs}, J_1^{\rs})\sim \mathcal{U}(\mathbb{N}_\rs^3)$, and $J_2^{\rs}\sim\mathcal{U}(\mathbb{N}_S\setminus\{J_1^\rs\})$
if $\{I_1^\rs\neq I_2^\rs\}$ and $J_2^\rs=J_1^\rs$ if $\{I_1^\rs= I_2^\rs\}$.\par  
Next, for $\pi\sim \mathcal{U}(\mathbb{S}_n)$ independent of $(I_1^{\rs}, I_2^{\rs}, J_1^{\rs}, J_2^{\rs})$, define 
\begin{align}
	I_3^\rs&\equiv \pi_\rs^{-1}(J_1^\rs),\quad I_4^\rs\equiv \pi_\rs^{-1}(J_2^\rs),\\
	J_3^\rs&\equiv \pi_\rs(I_1^\rs),\quad J_4^\rs\equiv \pi_\rs(I_2^\rs).
\end{align}
Then, 
$I_1^\rs=I_2^\rs$ if and only if $I_3^\rs=I_4^\rs$. 
For every quadruple $\mathbbm{i}_s=(i_1^s, i_2^s, i_3^s, i_4^s)\in\mathbb{N}_s^4$ satisfying 
$i_1^s=i_2^s$ if and only if $i_3^s=i_4^s$, let us define a permutation of $\mathbb{N}_s$:
\begin{equation}
	\tau_{\mathbbm{i}_s}(\alpha) 
	\equiv\begin{cases}
		i_4^s\quad\text{if}\ \alpha=i_1^s,\\
		i_3^s\quad\text{if}\  \alpha=i_2^s,\\
		\alpha\quad \text{if}\ \alpha\in \mathbb{N}_s\setminus \mathbbm{i}_s. 
	\end{cases}
\end{equation}

\begin{figure}
\begin{center}
\begin{equation}
\begin{tikzcd}
	{\bullet I_1^s} & {\bullet I_2^s} & {\bullet I_3^s} & {\bullet I_4^s} \\
	{\bullet J_1^s} & {\bullet J_2^s} & {\bullet J_3^s} & {\bullet J_4^s}
	\arrow["{{{                                                \pi_s^{-1}}}}"{marking, allow upside down, pos=0.2}, color={rgb,255:red,214;green,92;blue,92}, from=2-1, to=1-3]
	\arrow["{{{                               \pi_s^{-1}}}}"{marking, allow upside down, pos=0.8}, color={rgb,255:red,214;green,92;blue,92}, from=2-2, to=1-4]
	\arrow["{{{\pi_s}}}"{marking, allow upside down, pos=0.3}, from=1-1, to=2-3]
	\arrow["{{{\pi_s}}}"{marking, allow upside down, pos=0.8}, from=1-2, to=2-4]
	\arrow["{{{\tau_{\mathbb{I}_s}}}}", color={rgb,255:red,92;green,92;blue,214}, from=1-2, to=1-3]
	\arrow["{{{\tau_{\mathbb{I}_s}}}}", color={rgb,255:red,92;green,92;blue,214}, curve={height=-24pt}, from=1-1, to=1-4]
	\arrow["{{{t}}}", color={rgb,255:red,214;green,153;blue,92}, tail reversed, from=1-1, to=1-2]
\end{tikzcd}
\end{equation}
		\caption{Mappings of the indices when $\mathcal{S}=s$.}
		\label{fig: mapping}
\end{center}
\end{figure}
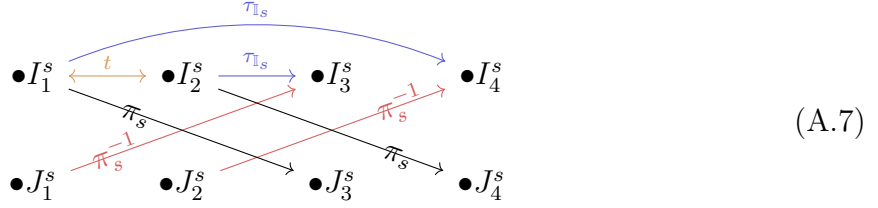
Furthermore, let $\mathbb{I}_s\equiv \{I_1^s, I_2^s, I_3^s, I_4^s\}, s=1,\dots,{S}$ and define 
\begin{align*}
	\pi^{\dagger}
	&\equiv \pi\circ \tau_{\mathbb{I}_\rs},\\
	\pi^{*}
	&\equiv \pi^{\dagger}\circ {t}({I}_1^\rs, I_2^\rs),
\end{align*}
where ${t}({I}_1^\rs, I_2^\rs)$ is the transposition of ${I}_1^\rs$ and $I_2^\rs$. A diagram illustrating the various mappings is depicted in Figure \ref{fig: mapping} for $\{\rs=s\}$ and $\bm{a}_{n}$ given. 
Then, following \cite{Bolthausen1984} and \cite{Schneller1988}, it can be shown that conditional on $\bm{a}_{n}$ and $\{\rs=s\}$,
\begin{enumerate}[label=(\roman*)]
	\item\label{property1} $\pi^\dagger_s(I_1^s)=J_2^s$,\ $\pi^\dagger_s(I_2^s)=J_1^s$,\ $\pi_s^{*}(I_1^s)=J_1^s$,\ $\pi_s^{*}(I_2^s)=J_2^s$;
	\item\label{property2}  $\pi_s, \pi_s^\dagger,\pi_s^{*}$ are independent of $\mathbb{I}_s$ and have identical $\mathcal{U}(\mathbb{G}_s)$ distribution;
	\item\label{property3}  $\pi_s^\dagger$ and $\{I_1^s, J_1^s\}$ are independent;
	\item\label{property4}  $(I_l^s, \tilde{\pi}_s(I_l^s))\sim\mathcal{U}(\mathbb{N}_s^2)$ for all $l\in\{1,2,3,4\}$ and $\tilde{\pi}_s\in\{\pi_s, \pi_s^\dagger,\pi_s^{*}\}$. 
\end{enumerate}
Finally, we let $\pi^\dagger\equiv (\pi_1^\dagger,\dots,\pi_{S}^\dagger)$ and $\pi^{*}\equiv (\pi_1^{*},\dots,\pi_{S}^{*})$.

\paragraph*{Step 2: Asymptotic normality}
By the mean value expansion, 
\begin{align}\label{eq: MVE}
	f(T^{\pi^{*}})
	&=f(T^{\pi^\dagger})+(T^{\pi^{*}}-T^{\pi^\dagger})\int_{0}^{1}f'\left(T^{\pi^\dagger}+t(T^{\pi^{*}}-T^{\pi^\dagger})\right)dt.
\end{align}
Rewrite \eqref{eq: MVE} as 
\begin{align}
	f(T^{\pi^{*}})
	&=f(T^{\pi^\dagger})+(T^{\pi^{*}}-T^{\pi^\dagger})f'(T^\pi)\notag\\
	&\quad +(T^{\pi^{*}}-T^{\pi^\dagger})\int_{0}^{1}\left(f'(T^{\pi^\dagger}+t(T^{\pi^{*}}-T^{\pi^\dagger}))-f'(T^\pi)\right)dt.\label{eq: mv expansion}
\end{align}
Next we remark that
\begin{align}
	n\E\left[\sigma_n^{-1}a_{I_1^{\rs}J_1^{\rs}}^{\rs}f'(T^{\pi^{*}})\big\vert \bm{a}_{n} \right]
	&=n\E\left[\E\left[\sigma_n^{-1}a_{I_1^{\rs}\pi^{*}(I_1^{\rs})}^{\rs}f'(T^{\pi^{*}})\vert \bm{a}_{n}, \rs,  \pi^{*}\right]\big\vert \bm{a}_{n}\right]\notag\\
	&=n\E\left[\sigma_n^{-1}\sum_{s=1}^{{S}}\E[a^{\rs}_{I_1^{\rs}\pi^{*}(I_1^{\rs})}f'(T^{\pi^{*}})\vert \bm{a}_{n}, {\rs}=s,  \pi^{*}]p_s\big\vert \bm{a}_{n}\right]\notag\\
	&=\E\left[\sigma_n^{-1}\sum_{s=1}^{{S}}\sum_{i=1}^{n_s}a_{i\pi^{*}(i)}^sf'(T^{\pi^{*}})\big\vert \bm{a}_{n}\right]\notag\\
	&=\E_{\pi^{*}}\left[T^{\pi^{*}}f'(T^{\pi^{*}})\right],\label{eq: ET}
\end{align}
where the first and second equalities are by the iterated expectations, the third equality holds because $\E[a^{\rs}_{I_1^{\rs}\pi^{*}(I_1^{\rs})}f'(T^{\pi^{*}})\vert \bm{a}_{n}, {\rs}=s, \pi^{*}]p_s
=n^{-1}\sum_{i=1}^{n_s}a_{i\pi^{*}(i)}^sf'(T^{\pi^{*}})$ by the property \ref{property2} above and the fact that $I_1^s\sim\mathcal{U}(\mathbb{N}_s)$, and the fourth is by the definition: $T^{\pi^{*}}\equiv \sigma_n^{-1}\sum_{s=1}^S\sum_{i=1}^{n_s}a_{i\pi^{*}(i)}^s$.  
From \eqref{eq: mv expansion} and \eqref{eq: ET}, we have 
\begin{align}
	\E_{\pi^{*}}\left[T^{\pi^{*}}f'(T^{\pi^{*}})\right]
	&=n\E\left[\sigma_n^{-1}a_{I_1^{\rs}J_1^{\rs}}^{\rs}f'(T^{\pi^{*}})\vert \bm{a}_{n}\right]\notag\\
	&=n\E\left[\sigma_n^{-1}a_{I_1^{\rs}J_1^{\rs}}^{\rs}f'(T^{\pi^\dagger})\vert \bm{a}_{n}\right]
	+n\E\left[\sigma_n^{-1}a_{I_1^{\rs}J_1^{\rs}}^{\rs}(T^{\pi^{*}}-T^{\pi^\dagger})f'(T^\pi)\vert \bm{a}_{n}\right]\notag\\
	&\quad +n\E\left[\sigma_n^{-1}a_{I_1^{\rs}J_1^{\rs}}^{\rs}(T^{\pi^{*}}-T^{\pi^\dagger})\int_{0}^{1}\left(f'(T^{\pi^\dagger}+t(T^{\pi^{*}}-T^{\pi^\dagger}))-f'(T^\pi)\right)dt\big\vert \bm{a}_{n}\right]\notag\\
	&\equiv A_1+A_2+A_3.\label{eq: E expansion}
\end{align}
Hereafter, we will analyze each summand on the RHS of \eqref{eq: E expansion} in three substeps. We will then complete the proof in the fourth substep.

\paragraph*{Substep 1: Calculation of $A_1$ in \eqref{eq: E expansion}}
We have 
\begin{align}\label{eq: E expansion a}
	A_1
	&=n\E\left[\sigma_n^{-1}a_{I_1^{\rs}J_1^{\rs}}^{\rs}f'(T^{\pi^\dagger})\vert \bm{a}_{n}\right]\notag\\
	&=n\,\sigma_n^{-1}\sum_{s=1}^{S}\E\left[a_{I_1^{\rs}J_1^{\rs}}^\rs f'(T^{\pi^\dagger})
	\vert \bm{a}_{n},\mathcal{S}=s\right]p_s\notag\\
	&=n\,\sigma_n^{-1}\sum_{s=1}^{S}\E\left[a_{I_1^{s}J_1^{s}}^s\vert \bm{a}_{n}\right]\E\left[f'(T^{\pi^\dagger})\vert \bm{a}_{n},\rs=s\right]p_s\notag\\
	&=0\ \ \text{a.s.},
\end{align}
where the first equality is by the definition, the second equality holds by the iterated expectations, the third holds by the independence of $\pi_s^\dagger$ and $\{I_1^s, J_1^s\}$ (the property \ref{property3}), and the last holds because $\E\left[a_{I_1^{s}J_1^{s}}^s\vert \bm{a}_{n}\right]=n_s^{-2}\sum_{i,j=1}^{n_s}a_{ij}^s=0\ \text{a.s.},$ which follows from the fact that 
$(I_1^s, J_1^s)\sim \mathcal{U}(\mathbb{N}_s^2)$ and Condition \ref{cond: centered} of the theorem. 

\paragraph*{Substep 2: Calculation of $A_2$ in \eqref{eq: E expansion}}
By the definition of $\pi^{\dagger}$, $\pi^{*}$ and the property \ref{property1},
\begin{align}
	T^{\pi^*}-T^{\pi^\dagger}
	&=\sigma_n^{-1}(a^{\rs}_{I_1^{\rs}\pi^\dagger(I_2^{\rs})} +a^{\rs}_{I_2^{\rs}\pi^\dagger(I_1^{\rs})} - a^{\rs}_{I_1^{\rs}\pi^\dagger(I_1^{\rs})} - a^{\rs}_{I_2^{\rs}\pi^\dagger(I_2^{\rs})})\notag\\
	&=\sigma_n^{-1}(a^{\rs}_{I_1^{\rs} J_1^{\rs}} +a^{\rs}_{I_2^{\rs} J_2^{\rs}} - a^{\rs}_{I_1^{\rs}J_2^{\rs}} - a^{\rs}_{I_2^\rs J_1^\rs}).\label{eq: first diff}
\end{align}
Since $\pi$ and $\{I_1^\rs, I_2^\rs, J_1^\rs, J_2^\rs\}$ are independent conditional on $\bm{a}_{n}$, 
it follows from \eqref{eq: first diff} that 
$\sigma_n^{-1}a_{I_1^{\rs}J_1^{\rs}}^\rs(T^{\pi^*}-T^{\pi^\dagger})$ and $T^\pi$ are also conditionally independent. Thus,
\begin{align}
	A_2
	&=n\E\left[\sigma_n^{-1}a_{I_1^{\rs}J_1^{\rs}}^\rs(T^{\pi^{*}}-T^{\pi^\dagger})f'(T^\pi)\big\vert \bm{a}_{n}\right]\notag\\
	&=n\E\left[\sigma_n^{-1}a_{I_1^{\rs}J_1^{\rs}}^\rs(T^{\pi^{*}}-T^{\pi^\dagger})\big\vert \bm{a}_{n}\right]\E\left[f'(T^\pi)\vert \bm{a}_{n}\right]\notag\\
	&=n\E\left[\sigma_n^{-2}a_{I_1^{\rs}J_1^{\rs}}^\rs(a^{\rs}_{I_1^\rs J_1^\rs} +a^{\rs}_{I_2^\rs J_2^\rs} - a^{\rs}_{I_1^\rs J_2^\rs} - a^{\rs}_{I_2^\rs J_1^\rs})\big\vert \bm{a}_{n}\right]\E\left[f'(T^\pi)\vert \bm{a}_{n}\right].\label{eq: EdTfprime}
\end{align}
By the iterated expectations,
\begin{align}\label{eq: 4 sums}
	&\E\left[\sigma_n^{-2}a_{I_1^{\rs}J_1^{\rs}}^\rs(a^{\rs}_{I_1^\rs J_1^\rs} +a^{\rs}_{I_2^\rs J_2^\rs} - a^{\rs}_{I_1^\rs J_2^\rs} - a^{\rs}_{I_2^\rs J_1^\rs})\vert \bm{a}_{n}\right]\notag\\
	&=\sum_{s=1}^{{S}}\E\left[\sigma_n^{-2}a_{I_1^{s}J_1^{s}}^s(a^{s}_{I_1^s J_1^s} +a^{s}_{I_2^s J_2^s} - a^{s}_{I_1^sJ_2^s} - a^{s}_{I_2^s J_1^s})\vert \bm{a}_{n}\right]p_s.
\end{align}
For the first summand of \eqref{eq: 4 sums}, we have 
\begin{align}
	\E[a_{I_1^{s}J_1^{s}}^{s2}\vert \bm{a}_{n}]
	&={n_s^{-2}}\sum_{i,j=1}^{n_s}a_{ij}^{s2}.\label{eq: 4 sums 1}
\end{align}
For the second summand of \eqref{eq: 4 sums}, 
\begin{align}
	\E\left[a_{I_1^{s}J_1^{s}}^sa^{s}_{I_2^s J_2^s}\vert \bm{a}_{n}\right]
	&=\E\left[a_{I_1^{s}J_1^{s}}^sa^{s}_{I_2^s J_2^s}\vert \bm{a}_{n}, J_1^s=J_2^s\right] P[J_1^s=J_2^s]
	+\E\left[a_{I_1^{s}J_1^{s}}^sa^{s}_{I_2^s J_2^s}\vert \bm{a}_{n}, J_1^s\neq J_2^s\right]P[J_1^s\neq J_2^s]\notag\\
	&=n_s^{-1}\E\left[a_{I_1^{s}J_1^{s}}^{s2}\vert \bm{a}_{n}\right]
	+n_s^{-1}(n_s-1)\E\left[a_{I_1^{s}J_1^{s}}^sa^{s}_{I_2^s J_2^s}\vert \bm{a}_{n}, J_1^s\neq J_2^s\right]\notag\\
	&=n_s^{-3}\sum_{i,j=1}^{n_s}a_{ij}^{s2}+n_s^{-1}(n_s-1)\E\left[a_{I_1^{s}J_1^{s}}^sa^{s}_{I_2^s J_2^s}\vert \bm{a}_{n}, J_1^s\neq J_2^s\right]\notag\\
	&=n_s^{-3}\sum_{i,j=1}^{n_s}a_{ij}^{s2}+n_s^{-1}(n_s-1)\frac{1}{n_s^2(n_s-1)^2}\sum_{i,j=1}^{n_s}a_{ij}^{s2}\notag\\
	&=\frac{1}{n_s^2(n_s-1)}\sum_{i,j=1}^{n_s}a_{ij}^{s2},\label{eq: 4 sums 2}
\end{align}
where the first equality is by the iterated expectations, the second equality holds because 
$a_{I_1^{s}J_1^{s}}^sa^{s}_{I_2^s J_2^s}=a_{I_1^{s}J_1^{s}}^{s2}$ if $J_1^s=J_2^s$, 
$P[J_1^s=J_2^s]=P[I_1^s=I_2^s]=n_s/n_s^2=n_s^{-1}$ and $P[J_1^s\neq J_2^s]=1-n_s^{-1}$, the third equality uses the fact that $(I_1^s, J_1^s)\sim \mathcal{U}(\mathbb{N}_s^2)$, the fifth is by algebra, and 
the fourth equality follows from the direct calculations below where 
Condition \ref{cond: centered} of the theorem is used in the third and fifth equalities:

\begin{align}
	\E\left[a_{I_1^{s}J_1^{s}}^sa^{s}_{I_2^s J_2^s}\vert \bm{a}_{n}, J_1^s\neq J_2^s\right]
	&=\frac{1}{n_s^2(n_s-1)^2}\sum_{i_1\neq i_2}\sum_{j_1\neq j_2}a_{i_1j_1}^sa_{i_2j_2}^s\notag\\
	&=\frac{1}{n_s^2(n_s-1)^2}\sum_{i_1\neq i_2}\sum_{j_1=1}^{n_s}a_{i_1j_1}^s\left(\sum_{j_2=1}^{n_s}a_{i_2j_2}^s-a_{i_2j_1}^s\right)\notag\\
	&=-\frac{1}{n_s^2(n_s-1)^2}\sum_{i_1\neq i_2}\sum_{j_1=1}^{n_s}a_{i_1j_1}^sa_{i_2j_1}^s\notag\\
	&=-\frac{1}{n_s^2(n_s-1)^2}\sum_{i_1=1}^{n_s}\left(\sum_{i_2=1}^{n_s}\sum_{j_1=1}^{n_s}a_{i_1j_1}^sa_{i_2j_1}^s-\sum_{j_1=1}^{n_s}a_{i_1j_1}^{s2}\right)\notag\\
	&=-\frac{1}{n_s^2(n_s-1)^2}\sum_{i_1=1}^{n_s}\sum_{i_2=1}^{n_s}\sum_{j_1=1}^{n_s}a_{i_1j_1}^sa_{i_2j_1}^s+\frac{1}{n_s^2(n_s-1)^2}\sum_{i,j=1}^{n_s}a_{ij}^{s2}\notag\\
	&=\frac{1}{n_s^2(n_s-1)^2}\sum_{i,j=1}^{n_s}a_{ij}^{s2}.\notag
\end{align}
For the third summand of \eqref{eq: 4 sums}, 
\begin{align}
	\E\left[a_{I_1^{s}J_1^{s}}^sa^{s}_{I_1^s J_2^s}\vert \bm{a}_{n}\right]
	&=\E\left[a_{I_1^{s}J_1^{s}}^sa^{s}_{I_1^s J_2^s}\vert \bm{a}_{n}, J_1^s=J_2^s\right] P[J_1^s=J_2^s]
	+\E\left[a_{I_1^{s}J_1^{s}}^sa^{s}_{I_1^s J_2^s}\vert \bm{a}_{n}, J_1^s\neq J_2^s\right]P[J_1^s\neq J_2^s]\notag\\
	&=n_s^{-1}\E\left[a_{I_1^{s}J_1^{s}}^{s2}\vert \bm{a}_{n}\right]
	+n_s^{-1}(n_s-1)\E\left[a_{I_1^{s}J_1^{s}}^sa^{s}_{I_1^s J_2^s}\vert \bm{a}_{n}, J_1^s\neq J_2^s\right]\notag\\
	&=n_s^{-3}\sum_{i,j=1}^{n_s}a_{ij}^{s2}-n_s^{-1}(n_s-1)\frac{1}{n_s^2(n_s-1)}\sum_{i=1}^{n_s}\sum_{j=1}^{n_s}a_{ij}^{s2}\notag\\
	&=0\ \ \text{a.s.},\label{eq: 4 sums 3}
\end{align}
where the second equality holds because $P[I_1^s=I_2^s]=P[J_1^s=J_2^s]=n_s^{-1}$ 
and $P[I_1^s\neq I_2^s]=P[J_1^s\neq J_2^s]=n_s^{-1}(n_s-1)$, and the third equality follows from 

\begin{align*}
	\E\left[a_{I_1^{s}J_1^{s}}^sa^{s}_{I_1^s J_2^s}\vert \bm{a}_{n}, J_1^s\neq J_2^s\right]
	&=n_s^{-1}\sum_{i=1}^{n_s}\E\left[a_{iJ_1^{s}}^sa^{s}_{iJ_2^s}\vert \bm{a}_{n}, J_1^s\neq J_2^s\right]\notag\\
	&=n_s^{-1}\sum_{i=1}^{n_s}\frac{1}{n_s(n_s-1)}\sum_{j_1\neq j_2}a_{ij_1}^sa_{ij_2}^s\notag\\
	&=n_s^{-1}\sum_{i=1}^{n_s}\frac{1}{n_s(n_s-1)}\sum_{j_1=1}^{n_s}a_{ij_1}^s\left(\sum_{j_2=1}^{n_s}a_{ij_2}^s-a_{ij_1}^s\right)\notag\\
	&=-\frac{1}{n_s^2(n_s-1)}\sum_{i=1}^{n_s}\sum_{j_1=1}^{n_s}a_{ij_1}^{s2},
\end{align*}
where the last equality uses Condition \ref{cond: centered} of the theorem. Similarly to 
\eqref{eq: 4 sums 3}, for the fourth summand of \eqref{eq: 4 sums}, we obtain
\begin{align}\label{eq: 4 sums 4}
	\E\left[a_{I_1^{s}J_1^{s}}^sa^{s}_{I_2^s J_1^s}\vert \bm{a}_{n}\right]
	=0\ \ \text{a.s.}
\end{align}
Using \eqref{eq: 4 sums 1}-\eqref{eq: 4 sums 4} in \eqref{eq: 4 sums}, we obtain 
\begin{align}
	&n\E\left[\sigma_n^{-2}a_{I_1^{\rs}J_1^{\rs}}^\rs(a^{\rs}_{I_1^\rs J_1^\rs} +a^{\rs}_{I_2^\rs J_2^\rs} - a^{\rs}_{I_1^\rs J_2^\rs} - a^{\rs}_{I_2^\rs J_1^\rs})\vert \bm{a}_{n}\right]\notag\\
	&=n\sigma_n^{-2}\sum_{s=1}^{{S}}\E\left[a_{I_1^{s}J_1^{s}}^s(a^{s}_{I_1^s J_1^s} +a^{s}_{I_2^s J_2^s} - a^{s}_{I_1^sJ_2^s} - a^{s}_{I_2^s J_1^s})\vert \bm{a}_{n}\right]p_s\notag\\
	&=n\sigma_n^{-2}\sum_{s=1}^{{S}}\left[\frac{1}{n_s^2}\sum_{i,j=1}^{n_s}a_{ij}^{s2}
	+\frac{1}{n_s^2(n_s-1)}\sum_{i,j=1}^{n_s}a_{ij}^{s2}\right]\frac{n_s}{n}\notag\\
	&=\sigma_n^{-2}\sum_{s=1}^{{S}}\frac{1}{n_s-1}\sum_{i,j=1}^{n_s}a_{ij}^{s2}\notag\\
	&=1\ \ \text{a.s.}\label{eq: E expansion b}
\end{align}
Using \eqref{eq: E expansion b} in \eqref{eq: EdTfprime}, we obtain 
\begin{align}
	A_2
	=\E_\pi\left[f'(T^\pi)\right].\label{eq: EdTfprime simple}
\end{align}
\paragraph*{Substep 1: Convergence of $A_3$ in \eqref{eq: E expansion}}
From the condition \eqref{eq: f prop2}, there exists a constant $M_0<\infty$ such that 
$|f'(x)|\leq M_0$ for all $x\in\mathbb{R}$, and 
since $f'(x)$ is uniformly continuous, 
for any $\eta>0$, there exists $\delta>0$ such that 
\begin{equation}\label{eq: f unicon}
\vert f'(x)-f'(y)\vert\leq \eta\quad\text{for all $x,y$ satisfying}\quad\vert x-y\vert \leq \delta.
\end{equation} 
By the convexity of $x\mapsto \vert x\vert$ and 
the inequality $\vert \int g(x)dx\vert \leq\int  \vert g(x)\vert dx$, we have 
\begin{align}
	A_3&
	=n\E\left[\sigma_n^{-1}a_{I_1^{\rs}J_1^{\rs}}^\rs(T^{\pi^{*}}-T^{\pi^\dagger})\int_{0}^{1}\left(f'(T^{\pi^\dagger}+t(T^{\pi^{*}}-T^{\pi^\dagger}))-f'(T^\pi)\right)dt\big\vert \bm{a}_{n}\right]\notag\\
	&\leq n\E\left[\sigma_n^{-1}\left\vert a_{I_1^{\rs}J_1^{\rs}}^\rs(T^{\pi^{*}}-T^{\pi^\dagger})\right\vert \left\vert\int_{0}^{1}\left(f'(T^{\pi^\dagger}+t(T^{\pi^{*}}-T^{\pi^\dagger}))-f'(T^\pi)\right)dt\right\vert  \big\vert \bm{a}_{n}\right]\notag\\
	&\leq n\E\bigg[\sigma_n^{-1}\left|a_{I_1^{\rs}J_1^{\rs}}^\rs(T^{\pi^{*}}-T^{\pi^\dagger})\right|\int_{0}^{1}|f'(T^{\pi^\dagger}+t(T^{\pi^{*}}-T^{\pi^\dagger}))-f'(T^\pi)|dt\notag\\
	&\qquad\qquad 1(|T^{\pi^{\dagger}}-T^\pi|+|T^{\pi^{*}}-T^{\pi^\dagger}|>\delta)\big\vert \bm{a}_{n}\bigg]\notag\\
	&\quad+n\E\bigg[\sigma_n^{-1}\left|a_{I_1^{\rs}J_1^{\rs}}^\rs(T^{\pi^{*}}-T^{\pi^\dagger})\right|\int_{0}^{1}|f'(T^{\pi^\dagger}+t(T^{\pi^{*}}-T^{\pi^\dagger}))-f'(T^\pi)|dt\notag\\
	&\qquad\qquad 1(|T^{\pi^{\dagger}}-T^\pi|+|T^{\pi^{*}}-T^{\pi^\dagger}|\leq \delta)\big\vert \bm{a}_{n}\bigg]\notag\\
	&\equiv {B}_1+B_2.\label{eq: E expansion c}	 	
\end{align} 
By the iterated expectations and the fact that $|f'(\cdot)|\leq M_0$, for $B_1$ in \eqref{eq: E expansion c}	 
\begin{align}
	B_1
	&=n\E\big[\sigma_n^{-1}\left|a_{I_1^{\rs}J_1^{\rs}}^\rs(T^{\pi^{*}}-T^{\pi^\dagger})\right|\int_{0}^{1}|f'(T^{\pi^\dagger}+t(T^{\pi^{*}}-T^{\pi^\dagger}))-f'(T^\pi)|dt\notag\\
	&\quad\quad 1\left(|T^{\pi^{\dagger}}-T^\pi|+|T^{\pi^{*}}-T^{\pi^\dagger}|>\delta\right)\big\vert {\bm{a}_{n}}\big]\notag\\
	&\leq 2M_0n\E\big[\sigma_n^{-1}\left|a_{I_1^{\rs}J_1^{\rs}}^\rs(T^{\pi^{*}}-T^{\pi^\dagger})\right|1\left(|T^{\pi^{\dagger}}-T^\pi|+|T^{\pi^{*}}-T^{\pi^\dagger}|>\delta\right)\big\vert \bm{a}_{n}\big],\notag\\
	&=2M_0\sum_{s=1}^{{S}}\E\big[\sigma_n^{-1}\left|a_{I_1^{s}J_1^{s}}^s(T^{\pi^{*}}-T^{\pi^\dagger})\right|1\left(|T^{\pi^{\dagger}}-T^\pi|+|(T^{\pi^{*}}-T^{\pi^\dagger}|>\delta\right)\big\vert \bm{a}_{n},\rs=s\big]n_s.\label{eq: B1 ubound}
\end{align}
From \eqref{eq: first diff} and the definitions of $\pi$ and $\pi^{\dagger}$, we have by the triangle inequality, conditional on $\bm{a}_{n}$ and $\{\rs=s\}$
\begin{align}
	|T^{\pi^{*}}-T^{\pi^{\dagger}}|+|T^{\pi^{\dagger}}-T^\pi|
	&\leq \sigma_n^{-1}\bigg(|a^{s}_{I_1^s J_1^s}| +|a^{s}_{I_2^s J_2^s}|+|a^{s}_{I_1^sJ_2^s}|+|a^{s}_{I_2^s J_1^s}|\notag\\
	&\quad+|a^{s}_{I_1^s \pi^\dagger(I_1^s)}| +|a^{s}_{I_2^s \pi^\dagger(I_2^s)}|+|a^{s}_{I_3^s \pi^\dagger(I_3^s)}|+|a^{s}_{I_4^s \pi^\dagger(I_4^s)}|\notag\\
	&\quad+|a^{s}_{I_1^s \pi(I_1^s)}| +|a^{s}_{I_2^s \pi(I_2^s)}|+|a^{s}_{I_3^s \pi(I_3^s)}|+|a^{s}_{I_4^s \pi(I_4^s)}|\bigg).\label{eq: 12 summands}	
\end{align}
If $|T^{\pi^{*}}-T^{\pi^{\dagger}}|+|T^{\pi^{\dagger}}-T^\pi|>\delta$, then at least 
one of the twelve summands in \eqref{eq: 12 summands} is greater than $\delta/12$ i.e. 
\begin{equation}\label{eq: aqr}
\sigma_n^{-1}|a_{QR}^s|>\delta/12\equiv \epsilon,
\end{equation}
 where $(Q, R)\sim\mathcal{U}(\mathbb{N}_s^2)$ by the property \ref{property4}. Thus, conditional on $\bm{a}_{n}$ and $\{\rs=s\}$
\begin{equation}\label{eq: B11}
	1(|T^{\pi^{\dagger}}-T^\pi|+|T^{\pi^{*}}-T^{\pi^\dagger}|>\delta)\leq 1(\sigma_n^{-1}|a_{QR}^s|>\epsilon).
\end{equation}
From \eqref{eq: first diff}, we have by the triangle inequality, conditional on $\bm{a}_{n}$ and $\{\rs=s\}$
\begin{align}\label{eq: B12}
	|T^{\pi^{\dagger}}-T^{\pi^{*}}|
	&\leq \sigma_n^{-1}\left(|a^{s}_{I_1^s J_1^s}| +|a^{s}_{I_2^s J_2^s}|+|a^{s}_{I_1^sJ_2^s}|+|a^{s}_{I_2^s J_1^s}|\right).
\end{align}
Using \eqref{eq: B11} and \eqref{eq: B12} in \eqref{eq: B1 ubound}, we obtain 
\begin{equation}
	B_1\leq 2M_0\sum_{s=1}^{S}n_s\E\left[\sigma_n^{-2}|a_{I_1^{s}J_1^{s}}^s|
	\left(|a^{s}_{I_1^s J_1^s}| +|a^{s}_{I_2^s J_2^s}|+|a^{s}_{I_1^sJ_2^s}|+|a^{s}_{I_2^s J_1^s}|\right)1(\sigma_n^{-1}|a_{QR}^s|>\epsilon)\vert \bm{a}_{n}\right].\label{eq: B1 ubound2}
\end{equation}
Next we will bound the terms $\E\left[\sigma_n^{-2}|a_{I_1^{s}J_1^{s}}^s||a^{s}_{GH}|1(\sigma_n^{-1}|a_{QR}^s|>\epsilon)\vert \bm{a}_{n}\right]$ in \eqref{eq: B1 ubound2},  
where $(G,H)\sim\mathcal{U}(\mathbb{N}_s^2)$. By a lengthy but straightforward algebra, 
\begin{align}
	&\E\left[\sigma_n^{-2}|a_{I_1^{s}J_1^{s}}^s||a^{s}_{GH}|1(\sigma_n^{-1}|a_{QR}^s|>\epsilon)\vert \bm{a}_{n}\right]\notag\\
	&=\E\left[\sigma_n^{-2}|a_{I_1^{s}J_1^{s}}^s||a^{s}_{GH}|1(\sigma_n^{-1}|a_{I_1^sJ_1^s}^s|>\epsilon)1(\sigma_n^{-1}|a_{GH}^s|>\epsilon)1(\sigma_n^{-1}|a_{QR}^s|>\epsilon)\vert \bm{a}_{n}\right]\notag\\
	&\quad+\E\left[\sigma_n^{-2}|a_{I_1^{s}J_1^{s}}^s||a^{s}_{GH}|1(\sigma_n^{-1}|a_{I_1^sJ_1^s}^s|\leq \epsilon)1(\sigma_n^{-1}|a_{GH}^s|> \epsilon)1(\sigma_n^{-1}|a_{QR}^s|>\epsilon)\vert \bm{a}_{n}\right]\notag\\
	&\quad+\E\left[\sigma_n^{-2}|a_{I_1^{s}J_1^{s}}^s||a^{s}_{GH}|1(\sigma_n^{-1}|a_{I_1^sJ_1^s}^s|>\epsilon)1(\sigma_n^{-1}|a_{GH}^s|\leq \epsilon)1(\sigma_n^{-1}|a_{QR}^s|>\epsilon)\vert \bm{a}_{n}\right]\notag\\
	&\quad+\E\left[\sigma_n^{-2}|a_{I_1^{s}J_1^{s}}^s||a^{s}_{GH}|1(\sigma_n^{-1}|a_{I_1^sJ_1^s}^s|\leq \epsilon)1(\sigma_n^{-2}|a_{GH}^s|\leq \epsilon)1(\sigma_n^{-1}|a_{QR}^s|>\epsilon)\vert \bm{a}_{n}\right]\notag\\
	&\leq \E\left[\sigma_n^{-2}|a_{I_1^{s}J_1^{s}}^s||a^{s}_{GH}|1(\sigma_n^{-1}|a_{I_1^sJ_1^s}^s|>\epsilon)1(\sigma_n^{-1}|a_{GH}^s|>\epsilon)\vert \bm{a}_{n}\right]\notag\\
	&\quad+\E\left[\sigma_n^{-2}|a^{s}_{GH}|\;\epsilon\; 1(\sigma_n^{-1}|a_{GH}^s|> \epsilon)1(\sigma_n^{-1}|a_{QR}^s|>\epsilon)\vert \bm{a}_{n}\right]\notag\\
	&\quad+\E\left[\sigma_n^{-2}|a_{I_1^{s}J_1^{s}}^s|1(\sigma_n^{-1}|a_{I_1^sJ_1^s}^s|>\epsilon)\; \epsilon\; 1(\sigma_n^{-1}|a_{QR}^s|>\epsilon)\vert \bm{a}_{n}\right]\notag\\
	&\quad+\E\left[\sigma_n^{-2}\epsilon^2 1(\sigma_n^{-1}|a_{QR}^s|>\epsilon)\vert \bm{a}_{n}\right]\notag\\
	&\leq 
	\E\left[\sigma_n^{-2}|a_{I_1^{s}J_1^{s}}^s||a^{s}_{GH}|1(\sigma_n^{-1}|a_{I_1^sJ_1^s}^s|>\epsilon)1(\sigma_n^{-1}|a_{GH}^s|>\epsilon)\vert \bm{a}_{n}\right]\notag\\
	&\quad+\E\left[\sigma_n^{-2}|a^{s}_{GH}||a_{QR}^s| 1(\sigma_n^{-1}|a_{GH}^s|> \epsilon)1(\sigma_n^{-1}|a_{QR}^s|>\epsilon)\vert \bm{a}_{n}\right]\notag\\
	&\quad+\E\left[\sigma_n^{-2}|a_{I_1^{s}J_1^{s}}^s||a^{s}_{QR}|1(\sigma_n^{-1}|a_{I_1^sJ_1^s}^s|>\epsilon)1(\sigma_n^{-1}|a_{QR}^s|>\epsilon)\vert \bm{a}_{n}\right]\notag\\
	&\quad+\E\left[\sigma_n^{-2}a_{QR}^{s2} 1(\sigma_n^{-1}|a_{QR}^s|>\epsilon)\vert \bm{a}_{n}\right]\notag\\
	&\leq 4 \E\left[\sigma_n^{-2}a_{I_1^{s}J_1^{s}}^{s2}1(\sigma_n^{-1}|a_{I_1^sJ_1^s}^s|>\epsilon)\vert \bm{a}_{n}\right]\notag\\
	&=4 \sigma_n^{-2}n_s^{-2}\sum_{i,j=1}^{n_s}a_{ij}^{s2} 1(\sigma_n^{-1}|a_{ij}^s|>\epsilon),\label{eq: B1 ubound3}
\end{align}
where the first equality holds by adding the four cases
$1(\sigma_n^{-1}|a_{I_1^sJ_1^s}^s|\gtreqless\epsilon)1(\sigma_n^{-1}|a_{GH}^s|\gtreqless\epsilon)$, 
 the first inequality uses $1(\sigma_n^{-1}|a_{QR}^s|>\epsilon)\leq 1$ and the inequalities in the indicator function $1(\cdot\leq \epsilon)$, the second inequality uses \eqref{eq: aqr}, the third inequality is by Cauchy-Schwarz inequality, and the last equality holds since $(I_1^s,J_1^s)\sim\mathcal{U}(\mathbb{N}_s^2)$.
Thus, using \eqref{eq: B1 ubound3} and Condition \ref{cond: Lindeberg} of the theorem in \eqref{eq: B1 ubound2} gives 
\begin{equation}\label{eq: B1 con}
	B_1\leq 32M_0\sigma_n^{-2}\sum_{s=1}^{{S}}n_s^{-1}\sum_{i,j=1}^{n_s}a_{ij}^{s2} 1(\sigma_n^{-1}|a_{ij}^s|>\epsilon)\conas 0.
\end{equation}
Next, we will bound $B_2$ as follows: 
\begin{align}
	B_2
	&=n\E\bigg[\sigma_n^{-1}\left|a_{I_1^{\rs}J_1^{\rs}}^\rs(T^{\pi^{*}}-T^{\pi^\dagger})\right|\int_{0}^{1}|f'(T^{\pi^\dagger}+t(T^{\pi^{*}}-T^{\pi^\dagger}))-f'(T^\pi)|dt\notag\\
	&\qquad 1(|T^{\pi^{\dagger}}-T^\pi|+|T^{\pi^{*}}-T^{\pi^\dagger}|\leq \delta)\big\vert \bm{a}_{n}\bigg]\notag\\
	&\leq \eta\, 
	n\E\left[\sigma_n^{-1}\left|a_{I_1^{\rs}J_1^{\rs}}^\rs(T^{\pi^{*}}-T^{\pi^\dagger})\right|1(|T^{\pi^{\dagger}}-T^\pi|+|T^{\pi^{*}}-T^{\pi^\dagger}|\leq \delta)\big\vert \bm{a}_{n}\right]\notag\\
	&\leq \eta\, 
	n\E\left[\sigma_n^{-1}\left|a_{I_1^{\rs}J_1^{\rs}}^\rs(T^{\pi^{*}}-T^{\pi^\dagger})\right|\big\vert \bm{a}_{n}\right]\notag\\
	&= \eta\, 
\sum_{s=1}^{{S}}n_s \E\left[\sigma_n^{-1}\left|a_{I_1^{\rs}J_1^{\rs}}^\rs(T^{\pi^{*}}-T^{\pi^\dagger})\right|\big\vert \bm{a}_{n}, \rs=s\right]\notag\\
	&\leq \eta\sum_{s=1}^{S}n_s \E\left[\sigma_n^{-2}|a_{I_1^{s}J_1^{s}}^s|(|a^{s}_{I_1^s J_1^s}| +|a^{s}_{I_2^s J_2^s}|+|a^{s}_{I_1^sJ_2^s}|+|a^{s}_{I_2^s J_1^s}|)\vert \bm{a}_{n}\right]\notag\\
	&\leq 4\eta \sum_{s=1}^{S}n_s \E\left[\sigma_n^{-2}a_{I_1^{s}J_1^{s}}^{s2}\vert \bm{a}_{n}\right]\notag\\
	&=4\eta\, \sigma_n^{-2}\sum_{s=1}^{S}n_s^{-1}\sum_{i,j=1}^{n_s}a_{ij}^{s2}\notag\\
	&\leq 4\eta\ \ \text{a.s.},
\end{align}
where the first equality is the definition of $B_2$, the first inequality is by \eqref{eq: f unicon}, the second inequality is by the inequality $1(\cdot)\leq 1$, the second equality is by the iterated expectations, the third inequality is by \eqref{eq: B12}, the fourth inequality is by Cauchy-Schwarz inequality, and the last holds by the definition of $\sigma_n^2$. Since $\eta$ is arbitrary,
\begin{equation}\label{eq: B2 con}
	B_2\conas 0.
\end{equation}  
Therefore, using \eqref{eq: B1 con} and \eqref{eq: B2 con} in \eqref{eq: E expansion c} gives 
\begin{equation}\label{eq: A3 con}
	A_3\conas 0.
\end{equation}
\paragraph*{Substep 4: Completing the proof}
To complete the proof of asymptotic normality, using  \eqref{eq: E expansion a}, \eqref{eq: EdTfprime simple} and \eqref{eq: A3 con} in \eqref{eq: E expansion}, we obtain $\E_{\pi^{*}}\left[T^{\pi^{*}}f(T^{\pi^{*}})\right]\conas \E_{\pi}\left[f'(T^{\pi})\right]$ as 
$n\to \infty$. Then, since $T^{\pi}$ and $T^{\pi^{*}}$ have identical distributions conditional on ${\bm{a}_{n}}$, it follows that 
\begin{align}
	\E_\pi\left[T^{\pi}f(T^{\pi})\right]\conas \E_{\pi}\left[f'(T^{\pi})\right].\label{eq: ETfprime con}
\end{align}
Finally, from \eqref{eq: f prop1} and \eqref{eq: ETfprime con}, we obtain $\E_\pi[h(T^\pi)]\conas \E[h(Z)]$ as 
$n\to \infty$, which, combined with the Portmanteau theorem, implies that 
$P^\pi(T^\pi\leq t)\conas \Phi(t)$ as required.\par

\section{Proof of Corollary \ref{cor: vector bc}}
Let $a_{ij}^s=n^{-1/2}t'\Sigma_n^{-1/2}\tilde{b}_{si}\tilde{c}_{sj}$ in Theorem \ref{HoeffdingCLT}, where $\tilde{b}_{si}\equiv b_{si}-\bar{b}_s$, $\bar{b}_s\equiv n_{s}^{-1}\sum_{i=1}^{n_s}b_{si}$, $b_{si}\in\mathbb{R}^k$ with $k\geq 1$, $t\in\mathbb{R}^k$ with $t\neq 0$,  $\tilde{c}_{sj}\equiv c_{sj}-\bar{c}_s$ and $\bar{c}_s\equiv n_{s}^{-1}\sum_{j=1}^{n_s}c_{sj}$. Then, for $\sigma_n^2$ in Theorem \ref{HoeffdingCLT}\ref{cond: Lindeberg}, we have 
\begin{align}
	\sigma_n^2
	&\equiv n^{-1}\sum_{s=1}^{{S}} \frac{1}{n_s-1}\left(\sum_{i=1}^{n_s}(t'\Sigma_n^{-1/2}\tilde{b}_{si})^{2}\right)
	\left(\sum_{i=1}^{n_s}\tilde{c}_{si}^{2}\right)=\Vert t\Vert^2>0,\label{eq: Vector bc sigma}
\end{align}
and the Lindeberg condition in Theorem \ref{HoeffdingCLT}\ref{cond: Lindeberg} becomes 
\begin{equation}\label{cond: Lindeberg prod}
	\sigma_n^{-2}\sum_{s=1}^{S} {n_s}^{-1}\left(n^{-1}\sum_{i,j=1}^{n_s}(t'\Sigma_n^{-1/2}\tilde{b}_{si})^{2}\tilde{c}_{sj}^{2}1\left(\sigma_n^{-1}n^{-1/2}|t'\Sigma_n^{-1/2}\tilde{b}_{si}\tilde{c}_{sj}|>\epsilon\right)\right)\conas 0.
\end{equation}
Under \eqref{cond: Lindeberg prod}, Theorem \ref{HoeffdingCLT} gives $
	n^{-1/2}\sum_{s=1}^{{S}}\sum_{i=1}^{n_s}t'\Sigma_n^{-1/2}\tilde{b}_{si}\tilde{c}_{s\pi(i)}\cond \Norm{0,1}.$ The latter combined with the Cram{\'e}r-Wold device yields \eqref{eq: Vector AN}. Below, we will verify \eqref{cond: Lindeberg prod} under Assumptions \ref{cond: Vector bc1}  and \ref{cond: Vector bc2} of the corollary. To this end, we first remark that 
\begin{equation}\label{eq: eval ineq}
\Vert t'\Sigma_n^{-1/2}\Vert^2\leq \Vert t\Vert^2(\lambda_{\min}(\Sigma_n))^{-1}
<\Vert t\Vert^2\lambda^{-1}.
\end{equation} 
\paragraph*{\eqref{cond: Lindeberg prod} via Assumption \ref{cond: Vector bc1} of Corollary \ref{cor: vector bc}} From the condition in \eqref{cond: Lyap}, \eqref{eq: Vector bc sigma}, \eqref{eq: eval ineq} and Cauchy-Schwarz inequality, 
	a sufficient condition for \eqref{cond: Lindeberg prod} is that for some $\delta>0$
	\begin{equation}\label{cond: Lyap bc}
		\frac{1}{\sigma_n^{2+\delta/2}}\sum_{s=1}^{S}\frac{1}{n_s}\sum_{i,j=1}^{n_s}|a_{ij}^{s}|^{2+\delta/2}
		\leq \frac{1}{\lambda^{1+\delta/4}n^{1+\delta/4}}		
		\sum_{s=1}^{{S}} \frac{1}{n_s}\left(\sum_{i=1}^{n_s}\Vert\tilde{b}_{si}\Vert^{2+\delta/2}\right)\left(\sum_{i=1}^{n_s}|\tilde{c}_{si}|^{2+\delta/2}\right)\conas 0.	
	\end{equation} 
From the $c_r$ and Jensen's inequalities, 
\begin{align}
&\Vert\tilde{b}_{si}\Vert^{4+\delta}\leq  2^{3+\delta}(\Vert{b}_{si}\Vert^{4+\delta}+n_s^{-1}\sum_{i=1}^{n_s}\Vert b_{si}\Vert^{4+\delta}),\label{eq: cr1}\\
&\vert\tilde{c}_{si}\vert^{4+\delta}\leq  2^{3+\delta}(\vert{c}_{si}\vert^{4+\delta}+n_s^{-1}\sum_{i=1}^{n_s}\vert c_{si}\vert^{4+\delta}).\label{eq: cr2}
\end{align}
For the RHS of \eqref{cond: Lyap bc},
	\begin{align}
		&n^{-1-\delta/4}\sum_{s=1}^{{S}} \left(\frac{1}{n_s}\sum_{i=1}^{n_s}\Vert\tilde{b}_{si}\Vert^{2+\delta/2}\right)
		\left(\frac{1}{n_s}\sum_{i=1}^{n_s}|\tilde{c}_{si}|^{2+\delta/2}\right)n_s\notag\\
		&\leq n^{-1-\delta/4}\left[\sum_{s=1}^{{S}} \left(\frac{1}{n_s}\sum_{i=1}^{n_s}\Vert \tilde{b}_{si}\Vert^{2+\delta/2}\right)^2n_s\right]^{1/2}
		\left[\sum_{s=1}^{{S}}\left(\frac{1}{n_s}\sum_{i=1}^{n_s}|\tilde{c}_{si}|^{2+\delta/2}\right)^2n_s\right]^{1/2}\notag\\
		&\leq n^{-1-\delta/4}\left[\sum_{s=1}^{{S}}\sum_{i=1}^{n_s}\Vert\tilde{b}_{si}\Vert^{4+\delta}\right]^{1/2}
		\left[\sum_{s=1}^{{S}}\sum_{i=1}^{n_s}|\tilde{c}_{si}|^{4+\delta}\right]^{1/2}\notag\\
		&\leq 2^{4+\delta}n^{-1-\delta/4}\left[\sum_{s=1}^{{S}}\sum_{i=1}^{n_s}\Vert{b}_{si}\Vert^{4+\delta}\right]^{1/2}
		\left[\sum_{s=1}^{{S}}\sum_{i=1}^{n_s}|{c}_{si}|^{4+\delta}\right]^{1/2}\notag\\
		&\conas 0, \label{cond: Lyap bc2}
	\end{align}
	where the first inequality is by Cauchy-Schwarz inequality, the second is by Jensen's inequality, and the third is by \eqref{eq: cr1} and \eqref{eq: cr2}, and the convergence follows from  \eqref{eq: Vector bc sigma} and Assumption \ref{cond: Vector bc1} of the corollary. Thus, \eqref{cond: Lyap bc} holds.
\paragraph*{\eqref{cond: Lindeberg prod} via Assumption \ref{cond: Vector bc2} of Corollary \ref{cor: vector bc}} 
From the inequality $\Vert a-b\Vert^2\leq 2(\Vert a\Vert^2+\Vert b\Vert^2)$ and 
Cauchy-Schwarz inequality, $n^{-1}\max_{s,i}\Vert \tilde{b}_{si}\Vert^2\leq 2n^{-1}\max_{s,i}\Vert {b}_{si}\Vert^2+2n^{-1}\max_{s,i}\Vert \bar{b}_{s}\Vert ^2\leq 4n^{-1}\max_{s,i}\Vert{b}_{si}\Vert^2$. Hence, 
\begin{align}
	n^{-1}\max_{s,i}\Vert \tilde{b}_{si}\Vert^2
	&\conas 0.\label{eq: maxcon}
\end{align}
For any $\epsilon, \epsilon_1>0$, choose $\epsilon_2$ such that $\frac{2^{4+\delta}M_0^2\epsilon_2^\delta}{\lambda^{1+\delta/2}\epsilon^\delta}
\leq \epsilon_1$. From \eqref{eq: maxcon}, there exists $n_0\in\mathbb{N}$ such that 
$n^{-1/2}\max_{s,i}\Vert\tilde{b}_{si}\Vert<\epsilon_2$ a.s. for all $n\geq n_0$. Then, 
\begin{align}\label{cond: Lyap2}
	&\frac{1}{\sigma_n^{2}n}\sum_{s=1}^{{S}} \frac{1}{n_s}\sum_{i,j=1}^{n_s}(t'\Sigma_n^{-1/2}\tilde{b}_{si})^2\tilde{c}_{sj}^{2}1\left(\sigma_n^{-1}|\tilde{c}_{sj}|n^{-1/2}|t'\Sigma_n^{-1/2}\tilde{b}_{si}|>\epsilon\right)\notag\\
		&\leq \frac{\Vert t\Vert^2}{\sigma_n^{2}\lambda n}\sum_{s=1}^{{S}} \frac{1}{n_s}\sum_{i,j=1}^{n_s}\Vert\tilde{b}_{si}\Vert^2\tilde{c}_{sj}^{2}1\left(\sigma_n^{-1}|\tilde{c}_{sj}|\Vert t\Vert\lambda^{-1/2}n^{-1/2}\Vert \max_{s,i} \Vert\tilde{b}_{si}\Vert >\epsilon\right)\notag\\
		&\leq \frac{1}{\lambda n}\sum_{s=1}^{{S}} \left(\sum_{i=1}^{n_s}\Vert\tilde{b}_{si}\Vert^2\right)\left(\frac{1}{n_s}\sum_{j=1}^{n_s}\tilde{c}_{sj}^{2}1\left(|\tilde{c}_{sj}| >\epsilon\lambda^{1/2}\epsilon_2^{-1}\right)\right)\notag\\
	&\leq \frac{\epsilon_2^\delta}{\lambda^{1+\delta/2}\epsilon^\delta}\left(n^{-1}\sum_{s=1}^{{S}} \sum_{i=1}^{n_s}\Vert\tilde{b}_{si}\Vert^2\right)\left(\max_{1\leq s\leq S}n_{s}^{-1}\sum_{i=1}^{n_s}|\tilde{c}_{si}|^{2+\delta}\right)\notag\\
	&\leq \frac{2^{4+\delta}M_0^2\epsilon_2^\delta}{\lambda^{1+\delta/2}\epsilon^\delta}\notag\\
	&\leq \epsilon_1,
\end{align} 
where the first inequality holds by Cauchy-Schwarz inequality, \eqref{eq: eval ineq} and taking the maximum in the indicator function, the second inequality holds by \eqref{eq: Vector bc sigma} and the fact that $n^{-1/2}\max_{s,i}\Vert\tilde{b}_{si}\Vert<\epsilon_2$ a.s., the third inequality holds by using the inequality in the indicator and taking the maximum over $s$, the fourth inequality uses 
Assumption \ref{cond: Vector bc2} of the corollary and the convexity of $x\mapsto \vert x\vert^{2+\delta}$, and 
the last holds by the choice of $\epsilon_2$. Thus, \eqref{cond: Lindeberg prod} holds. 	
\section{Justification of Remark \ref{remark: R2}}\label{sec:R2}
Recall that the Lindeberg condition in Theorem \ref{HoeffdingCLT}\ref{cond: Lindeberg} is 
\begin{equation}\label{cond: Lindeberg prod2}
	\sigma_n^{-2}\sum_{s=1}^{{S}} n_s^{-1}\left(\sum_{i,j=1}^{n_s}a_{ij}^{s2}1(\sigma_n^{-1}|a_{ij}^s|>\epsilon)\right)=\sigma_n^{-2}\sum_{s=1}^{{S}} n_s^{-1}\left(\sum_{i,j=1}^{n_s}\tilde{b}_{si}^{2}\tilde{c}_{sj}^{2}1(\sigma_n^{-1}|\tilde{b}_{si}\tilde{c}_{sj}|>\epsilon)\right)\conas 0.
\end{equation}
First, using \eqref{eq: ws Vs} and \eqref{eq: vareq}, we rewrite \eqref{eq: LYLind} as follows:
\begin{equation}\label{eq: FPLindeberg}
	\sum_{s=1}^{{S}} \frac{1}{n_s-1}\left(\sum_{i=1}^{n_s}\frac{p_s^2(n_s-n_{1s})(y_{si}-\bar{y}_s)^2}{n_{1s}n_s\sigma_n^2}1\left(\sigma_n^{-1}p_s|y_{si}-\bar{y}_s|>\epsilon n_{1s}\right)\right)\conas 0.
\end{equation}
We need to show that \eqref{eq: FPLindeberg} implies \eqref{cond: Lindeberg prod2} which can be seen from the 
following arguments:
\begin{align*}
	&\sigma_n^{-2}\sum_{s=1}^{{S}} \frac{1}{n_s}\left(\sum_{i,j=1}^{n_s}\tilde{b}_{si}^{2}\tilde{c}_{sj}^{2}1(\sigma_n^{-1}|\tilde{b}_{si}\tilde{c}_{sj}|>\epsilon)\right)\notag\\
	&=
	\sum_{s=1}^{{S}} \frac{1}{n_s}\bigg[\sum_{i=1}^{n_{1s}}\sum_{j=1}^{n_s}
	\left(	\frac{1}{n_{1s}}-\frac{1}{n_s}\right)^2\frac{p_s^2(y_{sj}-\bar{y}_s)^2}{\sigma_n^2}1\left(\sigma_n^{-1}\left(\frac{1}{n_{1s}}-\frac{1}{n_s}\right)p_s|y_{sj}-\bar{y}_s|>\epsilon\right)\notag\\
	&\quad\quad+\sum_{i=n_{1s}+1}^{n_{s}}\sum_{j=1}^{n_{s}}	\frac{1}{n_s^2}\frac{p_s^2(y_{sj}-\bar{y}_s)^2}{\sigma_n^2}1\left(\sigma_n^{-1}p_s|y_{sj}-\bar{y}_s|>\epsilon n_{s}\right)\bigg]\notag\\
	&=
	\sum_{s=1}^{{S}} \frac{1}{n_s}\bigg[\sum_{j=1}^{n_s}
	\frac{(n_{s}-n_{1s})^2}{n_{1s}n_s^2}\frac{p_s^2(y_{sj}-\bar{y}_s)^2}{\sigma_n^2}1\left(\sigma_n^{-1}p_s|y_{sj}-\bar{y}_s|>\epsilon\frac{n_{1s}n_s}{n_s-n_{1s}}\right)\notag\\
	&\quad\quad+\sum_{j=1}^{n_{s}}	\frac{n_s-n_{1s}}{n_s^2}\frac{p_s^2(y_{sj}-\bar{y}_s)^2}{\sigma_n^2}1\left(\sigma_n^{-1}p_s|y_{sj}-\bar{y}_s|>\epsilon n_{s}\right)\bigg]\notag\\
	&\leq 	\sum_{s=1}^{{S}} \frac{1}{n_s-1}\bigg[\sum_{j=1}^{n_s}
	\frac{(n_{s}-n_{1s})^2}{n_{1s}n_s^2}\frac{p_s^2(y_{sj}-\bar{y}_s)^2}{\sigma_n^2}1\left(\sigma_n^{-1}p_s|y_{sj}-\bar{y}_s|>\epsilon n_{1s}\right)\notag\\
	&\quad\quad+\sum_{j=1}^{n_{s}}	\frac{n_{s}-n_{1s}}{n_s^2}\frac{p_s^2(y_{sj}-\bar{y}_s)^2}{\sigma_n^2}1\left(\sigma_n^{-1}p_s|y_{sj}-\bar{y}_s|>\epsilon n_{1s}\right)\bigg]\notag\\	
	&=	\sum_{s=1}^{{S}} \frac{1}{n_s-1}\left[
	\sum_{j=1}^{n_s}
	\frac{n_{s}-n_{1s}}{n_{1s}n_s}\frac{p_s^2(y_{sj}-\bar{y}_s)^2}{\sigma_n^2}1\left(\sigma_n^{-1}p_s|y_{sj}-\bar{y}_s|>\epsilon n_{1s}\right)\right]\notag\\
	&\conas 0,
\end{align*}
where the first equality holds by \eqref{eq: LYb2} and breaking the sum $\sum_{i=1}^{n_{s}}$ into the two sums $\sum_{i=1}^{n_{1s}}$ and $\sum_{i=n_{1s}+1}^{n_s}$, the second equality holds because the $n_{1s}$ summands in $\sum_{i=1}^{n_{1s}}$ are identical and so are the $n_s-n_{1s}$ summands in $\sum_{i=n_{1s}+1}^{n_s}$, the first inequality holds using
$\frac{n_{1s}n_s}{n_s-n_{1s}}> n_{1s}$ and $n_{s}>n_{1s}$ in the indicator functions, and $n_s^{-1}<(n_s-1)^{-1}$ , the third equality is by simple algebra and, finally, the convergence holds
by \eqref{eq: FPLindeberg}.

\section{Proof of Proposition \ref{prop: SRIVa}}
	The strata of sizes $n_s=1$ do not contribute to the statistic in \eqref{def: T}. Since $(n-S)\to\infty$, without loss of generality, we may assume that $n_s\geq 2$, and the number of observations stratified into strata of size $n_s\geq 2$ tends to infinity.  Then, the proof follows by setting $b_{si}=q_{si}$ and $c_{si}=\rho_{si}$ in Corollary \ref{cor: vector bc}.  
\section{Proof of Proposition \ref{prop: AMstat}}
As in the proof of Proposition \ref{prop: SRIVa}, we can assume $n_s\geq 2$.
Let $b_{si}=Z_{si}$ and $c_{sj}=\varphi\left(\frac{j}{n_s+1}\right)$ in Corollary \ref{cor: vector bc}. 
From \eqref{eq: int ineq} and Assumption \ref{AM1} of Proposition \ref{prop: AMstat}, $n_s^{-1}\sum_{i=1}^{n_s}|\varphi(i/(n_s+1))|^{2+\delta}\leq \int_{0}^{1}|\varphi(x)|^{2+\delta}dx<(n_s+1)n_s^{-1}C_0\leq 2C_0.$ 
Under Assumption \ref{AM2} of Proposition \ref{prop: AMstat}, Lemma 13.1 of \cite{vanderVaart(1998)} implies that 
$\{R_{si}\}_{i=1}^{n_s}$ are uniformly distributed over $\{1,\dots, n_s\}$. We can then rewrite 
$A_{n}=n^{-1/2}\sum_{s=1}^S\sum_{i=1}^{n_s}\tilde{b}_{si}\tilde{c}_{s\pi(i)}$ with  $\pi\sim\mathcal{U}(\mathbb{S}_n)$. 
Therefore, Assumption \ref{cond: Vector bc2} of Corollary \ref{cor: vector bc} holds, and 
\begin{equation}
\Omega_{n}^{*-1/2}n^{1/2}A_{n}=\Omega_{n}^{*-1/2}n^{-1/2}\sum_{s=1}^S\sum_{i=1}^{n_s}\tilde{b}_{si}\tilde{c}_{s\pi(i)}\cond \Norm{0, I_k}.
\end{equation} 
Finally, the continuous mapping theorem gives $B_{n}^{*}\cond\chi^2_k$. 

\bibliographystyle{chicago}
\bibliography{SRCCLT}
	
\end{document}